# ON LIMITS OF SEQUENCES OF RESOLVENT KERNELS FOR SUBKERNELS

IGOR M. NOVITSKII

Abstract. In this paper, we approximate to continuous bi-Carleman kernels vanishing at infinity by sequences of their subkernels of Hilbert-Schmidt type and try to construct the resolvent kernels for these kernels as limits of sequences of the resolvent kernels for the approximating subkernels.

## 1. Introduction

In the general theory of integral equations of the second kind in $L^2 = L^2(\mathbb{R})$, i.e., equations of the form

$$f(s) - \lambda \int_{\mathbb{R}} \boldsymbol{T}(s,t) f(t)\, dt = g(s) \quad \text{for almost every } s \in \mathbb{R}, \tag{1.1}$$

it is customary to call an integral kernel $\boldsymbol{T}_{|\lambda}$ a *resolvent kernel for* $\boldsymbol{T}$ *at* $\lambda$ if the integral operator it induces on $L^2$ is the Fredholm resolvent $T_{|\lambda} := T(I - \lambda T)^{-1}$ of the integral operator $T$, which is induced on $L^2$ by the kernel $\boldsymbol{T}$. Once the resolvent kernel $\boldsymbol{T}_{|\lambda}$ has been constructed, one can express the $L^2$-solution $f$ to equation (1.1) in a direct and simple fashion as

$$f(s) = g(s) + \lambda \int_{\mathbb{R}} \boldsymbol{T}_{|\lambda}(s,t) g(t)\, dt \quad \text{for almost every } s \in \mathbb{R},$$

regardless of the particular choice of the function $g$ of $L^2$. Here it should be noted that, in general, the property of being an integral operator is not shared by Fredholm resolvents of integral operators, and there is even an example, given in [17] (see also [18, Section 5, Theorem 8]), of an integral operator having the property that at each non-zero regular value of the parameter $\lambda$, its Fredholm resolvent is not an integral operator. This phenomenon, however, can never occur for Carleman operators (the integral operators on $L^2$ whose kernels are Carleman, i.e., are square integrable in the second variable for almost all values of the first) due to the fact that the right-multiple by a bounded operator of a Carleman operator is again a Carleman operator (see Proposition 1 below). Therefore, in the case when the kernel $\boldsymbol{T}$ is Carleman and $\lambda$ is a regular value for $T$, the problem of solving equation (1.1) may be reduced to the problem of explicitly constructing in terms of $\boldsymbol{T}$ the resolvent kernel $\boldsymbol{T}_{|\lambda}$ which is a priori known to exist. A general formulation of this latter problem is given in Korotkov [16, Section 5, Problem 4] and is as follows:

($*$) *Let* $T$ *be a Carleman operator on* $L^2$*, and let* $\boldsymbol{T}$ *be the kernel of* $T$*. Develop methods of finding the set* $\Pi(T)$ *of all regular values of* $T$ *and the kernel* $\boldsymbol{T}_{|\lambda}$ *of the Fredholm resolvent* $T_{|\lambda} = T(I - \lambda T)^{-1}$ *at each point* $\lambda \in \Pi(T)$.

In order to approach this problem, we here confine our investigation to the case in which the kernel $\boldsymbol{T}\colon \mathbb{R}^2 \to \mathbb{C}$ of $T$ is *bi-Carleman* (i.e., square integrable in each variable separately for almost all values of the other) and is also a so-called $K^0$-*kernel* (i.e., it and its two Carleman functions $\boldsymbol{t}(s) = \overline{\boldsymbol{T}(s,\cdot)}$, $\boldsymbol{t'}(s) = \boldsymbol{T}(\cdot, s)\colon \mathbb{R} \to L^2$ are continuous and vanish at infinity; see Definition 1 below). Such conditions on $\boldsymbol{T}$ can always be achieved by means of a unitary equivalence transformation of $T$ if the $L^2$-adjoint to $T$ is also an integral operator (see Proposition 2 below). They therefore involve no loss of generality as far as the search of the $L^2$-solutions to integral equations of the form (1.1) is concerned (see Remarks 2 and 4 below).

A brief outline of our approach is as follows: we approximate to the $K^0$-kernels $\boldsymbol{T}$ by sequences of their subkernels $\boldsymbol{T}_n$ of Hilbert-Schmidt type (defined in Subsection 2.7 below) and introduce attempts to construct the resolvent kernels $\boldsymbol{T}_{|\lambda}$ for $\boldsymbol{T}$ as limits of sequences of the resolvent kernels $\boldsymbol{T}_{n|\lambda_n}$ for $\boldsymbol{T}_n$ where $\lambda_n \to \lambda$ as $n \to \infty$. General questions that are partially answered include the following: in what cases does $\{\boldsymbol{T}_{n|\lambda_n}\}$ converge and what is its limit? The convergence type used in most of the paper is with respect to sup-norms of appropriate spaces of continuous functions vanishing at infinity. The main results of the paper are contained in Theorems 3-5 and 7-9. They complement and improve the results of our previous paper [28].

## 2. Notation, definitions, and some background facts

2.1. **Spaces.** Throughout this paper, the symbols $\mathbb{C}$ and $\mathbb{N}$ refer to the complex plane and the set of all positive integers, respectively, $\mathbb{R}$ is the real line equipped with the Lebesgue measure, and $L^2 = L^2(\mathbb{R})$ is the complex Hilbert space of (equivalence classes of) measurable complex-valued functions on $\mathbb{R}$ equipped with the inner product $\langle f, g\rangle = \int f(s)\overline{g(s)}\,ds$ and the norm $\|f\| = \langle f, f\rangle^{\frac{1}{2}}$. (Integrals with no indicated domain, such as the above, are always to be extended over $\mathbb{R}$.) If $L \subset L^2$, we write $\overline{L}$ for the norm closure of $L$ in $L^2$, $L^{\perp}$ for the orthogonal complement of $L$ in $L^2$, and $\operatorname{Span}(L)$ for the norm closure of the set of all linear combinations of elements of $L$. Recall that a set $L$ in a normed space $Y$ is said to be *relatively compact* in $Y$ if each sequence of elements from $L$ contains a subsequence converging in the norm of $Y$. The sequence $\{f_n\}_{n=1}^{\infty}$ is also said to be *relatively compact* in $Y$ if the set of its values $\cup_{n=1}^{\infty} f_n$ is relatively compact in $Y$.

If $k$ is in $\mathbb{N}$ and $B$ is a Banach space with norm $\|\cdot\|_B$, let $C(\mathbb{R}^k, B)$ denote the Banach space, with the norm $\|f\|_{C(\mathbb{R}^k,B)} = \sup\limits_{x\in\mathbb{R}^k} \|f(x)\|_B$, of all continuous functions $f$ from $\mathbb{R}^k$ into $B$ that vanish at infinity, i.e., such that

$$\lim_{|x|\to\infty} \|f(x)\|_B = 0,$$

where $|\cdot|$ is the euclidian norm in $\mathbb{R}^k$. Given an equivalence class $f \in L^2$ containing a function of $C(\mathbb{R}, \mathbb{C})$, the symbol $[f]$ is used to mean that function.

A formulation of the Ascoli Theorem (see, e.g., [30, Theorem 1]) states:

**Theorem 1.** *The relative compactness of the set $L$ in $C(\mathbb{R}^k, B)$ is equivalent to the conjunction of two conditions:*

(A1) *the set $\bigcup\limits_{f\in L} f(x)$ is, for each fixed $x \in \mathbb{R}^k$, a relatively compact set in $B$,*

(A2) *$L$ is equicontinuous on $\mathbb{R}^k \cup \{\infty\}$, i.e., given $\varepsilon > 0$, there exists for every fixed $x \in \mathbb{R}^k$ a number $\delta_x(\varepsilon) > 0$ such that for $y$ in $\mathbb{R}^k$,*

$$\sup_{f\in L} \|f(x) - f(y)\|_B < \varepsilon \quad \text{whenever } |x - y| < \delta_x(\varepsilon), \tag{2.1}$$

*and there exists a number* $\delta_\infty(\varepsilon) > 0$ *such that for* $y$ *in* $\mathbb{R}^k$,

$$\sup_{f\in L} \|f(y)\|_B < \varepsilon \quad \textit{whenever } |y| > \delta_\infty(\varepsilon).$$

2.2. **Operators.** Let $\mathfrak{R}(L^2)$ denote the Banach algebra of all bounded linear operators from $L^2$ into itself; $\|\cdot\|$ will also denote the norm in $\mathfrak{R}(L^2)$. For an operator $A$ of $\mathfrak{R}(L^2)$, $A^*$ stands for the adjoint to $A$ with respect to $\langle\cdot,\cdot\rangle$, $\operatorname{Ran} A = \cup_{f\in L^2} Af$ for the range of $A$, and $\operatorname{Ker} A = \{f \in L^2 \mid Af = 0\}$ for the null-space of $A$. An operator $U \in \mathfrak{R}(L^2)$ is said to be *unitary* if $\operatorname{Ran} U = L^2$ and $\langle Uf, Ug\rangle = \langle f, g\rangle$ for all $f$, $g \in L^2$. An operator $A \in \mathfrak{R}(L^2)$ is said to be *invertible* if it has an inverse which is also in $\mathfrak{R}(L^2)$, i.e., if there is an operator $B \in \mathfrak{R}(L^2)$ for which $BA = AB = I$, where $I$ is the identity operator on $L^2$; $B$ is denoted by $A^{-1}$. An operator $P \in \mathfrak{R}(L^2)$ is called a *projection* in $L^2$ if $P^2 = P$. A projection $P$ in $L^2$ is said to be *orthogonal* if $P = P^*$, and let $\mathfrak{P}(L^2)$ denote the set of all orthogonal projections in $L^2$. An operator $T \in \mathfrak{R}(L^2)$ is said to be *self-adjoint* if $T^* = T$. An operator $T \in \mathfrak{R}(L^2)$ is said to be *compact* if it transforms every bounded set in $L^2$ into a relatively compact set in $L^2$. A (compact) operator $A \in \mathfrak{R}(L^2)$ is *nuclear* if $\sum_n |\langle Au_n, u_n\rangle| < \infty$ for any choice of an orthonormal basis $\{u_n\}$ of $L^2$.

An operator $T \in \mathfrak{R}(L^2)$ is said to be *asymptotically quasi-compact* if there is a sequence $\{S_n\}$ of compact operators such that $\|T^n - S_n\|^{1/n} \to 0$ as $n \to \infty$; if moreover, there is a positive integer $m$ such that $T^m$ is compact, we call $T$ simply a *quasi-compact* operator (both the terms are borrowed from [32, p. 19]).

2.3. **Fredholm resolvents.** Throughout this and the next subsection, $T$ denotes a bounded linear operator of $\mathfrak{R}(L^2)$. The set of *regular values* for $T$ (see [32, p. 28]), denoted by $\Pi(T)$, is the set of complex numbers $\lambda$ such that the operator $I - \lambda T$ is invertible, i.e., it has an inverse $R_\lambda(T) = (I - \lambda T)^{-1}$ in $\mathfrak{R}(L^2)$ that satisfies

$$(I - \lambda T)\, R_\lambda(T) = R_\lambda(T)\, (I - \lambda T) = I. \tag{2.2}$$

The operator

$$T_{|\lambda} := T R_\lambda(T) \ (= R_\lambda(T) T) \tag{2.3}$$

is then referred to as the *Fredholm resolvent* of $T$ at $\lambda$. Remark that if $\lambda$ is a regular value for $T$, then, for each fixed $g$ in $L^2$, the (unique) solution $f$ of $L^2$ to the second-kind equation $f - \lambda T f = g$ may be written as

$$f = g + \lambda T_{|\lambda} g$$

(follows from the formula

$$R_\lambda(T) = I + \lambda T_{|\lambda} \tag{2.4}$$

which is a rewrite of (2.2)). Recall that the inverse $R_\lambda(T)$ of $I - \lambda T$ as a function of $T$ also satisfies the following identity, often referred to as the *second resolvent equation* (see, e.g., [11, Theorem 5.16.1]): for $T$, $A \in \mathfrak{R}(L^2)$,

$$\begin{aligned} R_\lambda(T) - R_\lambda(A) &= \lambda R_\lambda(T)(T - A)R_\lambda(A) \\ &= \lambda R_\lambda(A)(T - A)R_\lambda(T) \quad \text{for every } \lambda \in \Pi(T) \cap \Pi(A). \end{aligned} \tag{2.5}$$

(A slightly modified version of it is

$$\begin{aligned} T_{|\lambda} - A_{|\lambda} &= (I + \lambda T_{|\lambda})(T - A)(I + \lambda A_{|\lambda}) \\ &= (I + \lambda A_{|\lambda})(T - A)(I + \lambda T_{|\lambda}) \quad \text{for every } \lambda \in \Pi(T) \cap \Pi(A), \end{aligned} \tag{2.6}$$

which involves the Fredholm resolvents.) It should also be mentioned that the map $R_\lambda(T)\colon \Pi(T) \to \mathfrak{R}(L^2)$ is continuous at every point $\lambda$ of the open set $\Pi(T)$, in the sense that

$$\|R_{\lambda_n}(T) - R_\lambda(T)\| \to 0 \tag{2.7}$$

when $\lambda_n \to \lambda$, $\lambda_n \in \Pi(T)$ (see, e.g., [13, Lemma 2 (XIII.4.3)]). Moreover, $R_\lambda(T)$ is given by an operator-norm convergent series ($T^0 = I$):

$$R_\lambda(T) = \sum_{n=0}^{\infty} \lambda^n T^n \quad \text{provided } |\lambda| < r(T) := \frac{1}{\lim\limits_{n\to\infty} \sqrt[n]{\|T^n\|}} \tag{2.8}$$

(see, e.g., [13, Theorem 1 (XIII.4.2)]).

Also note that if $T$ is self-adjoint and $\operatorname{Im}\lambda \neq 0$, then (see, e.g., [7, Lemma XII.2.2, p. 1192])

$$\|R_\lambda(T)\| \leqslant \frac{|\lambda|}{|\operatorname{Im}\lambda|}. \tag{2.9}$$

To simplify the formulae, we shall always write $R^*_\lambda(T)$ for the adjoint $(R_\lambda(T))^*$ to $R_\lambda(T)$.

Given a sequence $\{S_n\}_{n=1}^\infty$ of bounded operators on $L^2$, let $\nabla_{\mathfrak{b}}(\{S_n\})$ denote the set of all non-zero complex numbers $\zeta$ for which there exist positive constants $M(\zeta)$ and $N(\zeta)$ such that

$$\zeta \in \Pi(S_n) \text{ and } \left\|S_{n|\zeta}\right\| \leqslant M(\zeta) \quad \text{for } n > N(\zeta), \tag{2.10}$$

where, as in what follows, $S_{n|\zeta}$ stands for the Fredholm resolvent at $\zeta$ of the sequence term $S_n$, and let $\nabla_{\mathfrak{s}}(\{S_n\})$ denote the set of all non-zero complex numbers $\zeta$ ($\in \nabla_{\mathfrak{b}}(\{S_n\})$) for which the sequence $\left\{S_{n|\zeta}\right\}$ is convergent in the strong operator topology (i.e., the limit $\lim\limits_{n\to\infty} S_{n|\zeta} f$ exists in $L^2$ for every $f \in L^2$).

*Remark* 1. The sets $\nabla_{\mathfrak{b}}(\{S_n\})$ and $\nabla_{\mathfrak{s}}(\{S_n\})$ evidently remain unchanged if in their definition the Fredholm resolvents $S_{n|\zeta}$ are replaced by the operators $R_\zeta(S_n) = (I - \zeta S_n)^{-1} = I + \zeta S_{n|\zeta}$ (see (2.4)). So, if $\Delta_{\mathrm{b}}$ (resp., $\Delta_{\mathrm{s}}$) is the *region of boundedness* (resp., *strong convergence*) for the resolvents $\left\{(\zeta I - S_n)^{-1}\right\}$, which was introduced and studied in [14, Section VIII-1.1], then the sets $\nabla_{\mathfrak{b}}(\{S_n\})$ and $\Delta_{\mathrm{b}} \setminus \{0\}$ (resp., $\nabla_{\mathfrak{s}}(\{S_n\})$ and $\Delta_{\mathrm{s}} \setminus \{0\}$) are mapped onto each other by the mapping $\zeta \to \zeta^{-1}$. In the course of our investigations here, we keep this mapping in mind when referring to [14] for generalized strong convergence theory. In the next subsection, the mapping will also be of help to relate the spectrum with the characteristic set.

2.4. **Characteristic sets.** The *characteristic set* $\Lambda(T)$ for $T$ is defined to be the complementary set in $\mathbb{C}$ of $\Pi(T)$:

$$\Lambda(T) = \mathbb{C} \setminus \Pi(T) \tag{2.11}$$

(the name is adopted from [13, XIII.3.1]). Just as the spectrum $\sigma(T)$ of $T$ decomposes as the disjoint union $\sigma(T) = \sigma_p(T) \cup \sigma_c(T) \cup \sigma_r(T)$ of its point, continuous, and residual spectrums (see, e.g., [8, Definition XV.8.1]), the characteristic set $\Lambda(T)$ decomposes into the union

$$\Lambda(T) = \Lambda_p(T) \cup \Lambda_c(T) \cup \Lambda_r(T) \tag{2.12}$$

of the disjoint subsets

$$\begin{aligned}
\Lambda_p(T) &:= \{\lambda \in \Lambda(T) \mid \operatorname{Ker}\,(I - \lambda T) \neq \{0\}\} \\
&\quad\; \left(= \left\{\lambda \in \Lambda(T) \mid \tfrac{1}{\lambda} \in \sigma_p(T)\right\}\right), \\
\Lambda_c(T) &:= \left\{\lambda \in \Lambda(T) \mid \operatorname{Ker}\,(I - \lambda T) = \{0\} \text{ and } \overline{\operatorname{Ran}\,(I - \lambda T)} = L^2\right\} \\
&\quad\; \left(= \left\{\lambda \in \Lambda(T) \mid \tfrac{1}{\lambda} \in \sigma_c(T)\right\}\right), \\
\Lambda_r(T) &:= \left\{\lambda \in \Lambda(T) \mid \operatorname{Ker}\,(I - \lambda T) = \{0\} \text{ and } \overline{\operatorname{Ran}\,(I - \lambda T)} \neq L^2\right\} \\
&\quad\; \left(= \left\{\lambda \in \Lambda(T) \mid \tfrac{1}{\lambda} \in \sigma_r(T)\right\}\right).
\end{aligned} \tag{2.13}$$

Points of the set $\Lambda_p(T)$ are said to be *characteristic values* for $T$ (see again [13, XIII.3.1]).

Similarly as for the resolvent set and the respective parts of the spectrum in [37, Chapter 10, § 9] (or [12, Satz 5.13, S. 67]), it can be established that

$$(\Pi(T))^{\#} = \Pi(T^*), \quad (\Lambda_c(T))^{\#} = \Lambda_c(T^*), \quad (\Lambda_r(T))^{\#} \subset \Lambda_p(T^*), \tag{2.14}$$

where the superscript $^{\#}$ stands for the image under complex conjugation of a set in $\mathbb{C}$. Also notice that simultaneously

$$\Lambda_c(T) = \varnothing, \quad \Lambda_r(T) = \varnothing, \quad \Lambda(T) = \Lambda_p(T), \tag{2.15}$$

whenever $T$ is an asymptotically quasi-compact operator (see, e.g., [29, Proposition 4.1.5, p. 169]).

Let $T \in \mathfrak{R}(L^2)$ be a self-adjoint operator and $E_T(\cdot)\colon \Omega \to \mathfrak{P}(L^2)$ be its associated resolution of the identity, defined on the $\sigma$-algebra $\Omega$ of all Borel sets in $\mathbb{R}$ and supported on the spectrum $\sigma(T)$. Then the subdivisions of $\Lambda(T)$ displayed in (2.13) admit the following description in terms of $E_T(\cdot)$:

$$\begin{gathered} \Lambda_p(T) = \left\{\lambda \in \mathbb{R} \mid E_T\left(\left\{\tfrac{1}{\lambda}\right\}\right) \text{ is a non-zero projection in } L^2\right\}, \\ \Lambda_c(T) = \left\{\lambda \in \Lambda(T) \mid E_T\left(\left\{\tfrac{1}{\lambda}\right\}\right) \text{ is the zero projection in } L^2\right\}, \\ \Lambda_r(T) = \varnothing \end{gathered} \tag{2.16}$$

(see, e.g., [8, Section XV.8] or [9, Theorem 1, p. 56]). Moreover,

$$\{\lambda \in \Pi(T) \mid \operatorname{Im}\lambda = 0\} = \{\lambda \mid \frac{1}{\lambda} \in \varrho(T)\} \cup \{0\}, \tag{2.17}$$

where $\varrho(T)$ (the resolvent set of $T$) is the union of all open intervals $\omega$ in $\mathbb{R}$ at which $E_T(\omega)$ is the zero projection.

2.5. **Integral operators.** A linear operator $T\colon L^2 \to L^2$ is *integral* if there is a complex-valued measurable function $\boldsymbol{T}$ (*kernel*) on the Cartesian product $\mathbb{R}^2 = \mathbb{R} \times \mathbb{R}$ such that

$$(Tf)(s) = \int \boldsymbol{T}(s,t) f(t)\, dt \tag{2.18}$$

for every $f$ in $L^2$ and almost every $s$ in $\mathbb{R}$. Recall (see [10, Theorem 3.10]) that integral operators are bounded from $L^2$ into itself, and need not be compact.

A measurable function $\boldsymbol{T}\colon \mathbb{R}^2 \to \mathbb{C}$ is said to be *Carleman* (resp., *Hilbert-Schmidt*) kernel if $\int |\boldsymbol{T}(s,t)|^2\, dt < \infty$ for almost every fixed $s$ in $\mathbb{R}$ (resp., $\int\int |\boldsymbol{T}(s,t)|^2\, dt\, ds < \infty$). The nuclear operators of $\mathfrak{R}(L^2)$ are good examples of the integral operators whose kernels are Hilbert-Schmidt (see, e.g., [31, Theorem VI.23, pp. 210–211]).

To each Carleman kernel $\boldsymbol{T}$ there corresponds a *Carleman function* $\boldsymbol{t}\colon \mathbb{R} \to L^2$ defined by $\boldsymbol{t}(s) = \overline{\boldsymbol{T}(s,\cdot)}$ for all $s$ in $\mathbb{R}$ for which $\boldsymbol{T}(s,\cdot) \in L^2$. The Carleman kernel $\boldsymbol{T}$ is called *bi-Carleman* in case its conjugate transpose kernel $\boldsymbol{T'}$ ($\boldsymbol{T'}(s,t) = \overline{\boldsymbol{T}(t,s)}$) is also a Carleman kernel. Associated with the conjugate transpose $\boldsymbol{T'}$ of every bi-Carleman kernel $\boldsymbol{T}$ there is therefore an (also primed) Carleman function $\boldsymbol{t'}\colon \mathbb{R} \to L^2$ defined by $\boldsymbol{t'}(s) = \overline{\boldsymbol{T'}(s,\cdot)}\,(= \boldsymbol{T}(\cdot,s))$ for all $s$ in $\mathbb{R}$ for which $\boldsymbol{T'}(s,\cdot) \in L^2$. With each bi-Carleman kernel $\boldsymbol{T}$, we therefore associate the pair of Carleman functions $\boldsymbol{t}$, $\boldsymbol{t'}\colon \mathbb{R} \to L^2$, both defined, via $\boldsymbol{T}$, as above.

An integral operator whose kernel is Carleman (resp., bi-Carleman) is referred to as the *Carleman* (resp., *bi-Carleman*) operator. The integral operator $T$ is called *bi-integral* if its adjoint $T^*$ is also an integral operator; in that case if $\boldsymbol{T^*}$ is the kernel of $T^*$ then, in the above notation, $\boldsymbol{T^*}(s,t) = \boldsymbol{T'}(s,t)$ for almost all $(s,t) \in \mathbb{R}^2$ (see, e.g., [10, Theorem 7.5]). A bi-Carleman operator is always a bi-integral operator, but not conversely. Note also that not every bi-Carleman kernel is the kernel for an integral operator defined on all of $L^2$ as in (2.18).

*Remark* 2. The bi-integral operators, $T$, are generally involved in second-kind integral equations in $L^2$, of the form $f-\lambda Tf=g$, because their adjoint equations, which are of the form $u-\overline{\lambda}T^*u=v$, are typically required to be also integral and of the same kind. As our major concern in this paper is to try to devise some method to solve just such integral equations, in our further investigations we shall mostly be concerned with integral kernels which generate bi-integral operators of $\mathfrak{R}(L^2)$.

We conclude this subsection by recalling an important algebraic property of Carleman operators which will be exploited frequently throughout the text, a property which is the content of the following so-called "Right-Multipilication Lemma" (see [19], [15, Corollary IV.2.8], or [10, Theorem 11.6]):

**Proposition 1.** *Let $T$ be a Carleman operator, let $\boldsymbol{t}$ be the Carleman function associated with the inducing Carleman kernel of $T$, and let $A\in\mathfrak{R}(L^2)$ be arbitrary. Then the product operator $P=TA$ is also a Carleman operator, and the composition function*

$$\boldsymbol{p}(\cdot)=A^*(\boldsymbol{t}(\cdot))\colon \mathbb{R}\to L^2 \tag{2.19}$$

*is the Carleman function associated with its kernel.*

Before leaving this subsection, let us also remark that throughout the following we shall continue to employ the convention of referring to integral operators by (possibly tilded, subscripted, etc.) italic caps and to the corresponding kernels (resp., Carleman functions) by the same letter, but written in upper case (resp., lower case) bold-face type. Thus, e.g., if $\widetilde{T}_\alpha$ denotes, say, a bi-Carleman operator, then $\widetilde{\boldsymbol{T}}_\alpha$ and $\widetilde{\boldsymbol{t}}_\alpha$ and $\widetilde{\boldsymbol{t}}'_\alpha$ are to be used to denote respectively its kernel and two Carleman functions, and vice versa.

2.6. **$\boldsymbol{K}^0$-kernels.** Among all possible bi-Carleman kernels on $\mathbb{R}^2$, the following definition distinguishes a special type: those which together with their associated Carleman functions have the property of being continuous vanishing at infinity.

**Definition 1.** A bi-Carleman kernel $\boldsymbol{T}\colon \mathbb{R}^2\to\mathbb{C}$ is called a $K^0$-*kernel* if the following three conditions are satisfied:

(i) the function $\boldsymbol{T}$ is in $C(\mathbb{R}^2,\mathbb{C})$,
(ii) the Carleman function $\boldsymbol{t}$ associated with $\boldsymbol{T}$, $\boldsymbol{t}(s)=\overline{\boldsymbol{T}(s,\cdot)}$, is in $C(\mathbb{R},L^2)$,
(iii) the Carleman function $\boldsymbol{t'}$ associated with the conjugate transpose $\boldsymbol{T'}$ of $\boldsymbol{T}$, $\boldsymbol{t'}(s)=\overline{\boldsymbol{T'}(s,\cdot)}=\boldsymbol{T}(\cdot,s)$, is in $C(\mathbb{R},L^2)$.

Of course, this definition does not exclude the possibility of the existence of $K^0$-kernels which are associated with no integral operator belonging to $\mathfrak{R}(L^2)$. We also emphasize that the introduced kernels are special even among the set of all continuous bi-Carleman kernels on $\mathbb{R}^2$. The point of the following nontrivial example is to show that there are in $\mathfrak{R}(L^2)$ bi-Carleman operators that are induced by continuous kernels having both their Carleman functions continuous, but possessing none of the properties (i)-(iii).

*Example* 1. Let $\varepsilon$ be a fixed positive number and put $c_\varepsilon:=\frac{e^{-\varepsilon}}{\varepsilon}$. Define for $t\geqslant 0$ and $s\in\mathbb{R}$, the quantities $t_\varepsilon:=\max\{t-\varepsilon,0\}$ and $|s|_\varepsilon:=|s|+\varepsilon$, and the intervals $\mathbb{I}_s:=[|s|,|s|_\varepsilon]$, $\widehat{\mathbb{I}}_s:=(|s|_\varepsilon,+\infty)$, $\mathbb{I}_t^+:=[t_\varepsilon,t]$, $\mathbb{I}_t^-:=[-t,-t_\varepsilon]$ and $\mathbb{I}_t^\pm:=(\varepsilon-t,t-\varepsilon)$. Then define a bi-Carleman kernel $\boldsymbol{K}:\mathbb{R}^2\to[0,\infty)$ by writing

$$\boldsymbol{K}(s,t)=\begin{cases} e^{|s|-t} & \text{if } t\in\widehat{\mathbb{I}}_s,\\ -c_\varepsilon(|s|-t) & \text{if } t\in\mathbb{I}_s,\\ 0 & \text{otherwise.}\end{cases} \tag{2.20}$$

The inequalities between the extreme terms of the following two chains of relations imply by the Schur test (see, e.g., [10, Theorem 5.2]) that there exists in $\mathfrak{R}(L^2)$ an integral operator $K$ whose kernel is just the above-defined function $\boldsymbol{K}$:

$$\begin{aligned}\int \boldsymbol{K}(s,t)q(t)\,dt &= -c_\varepsilon \int_{\mathbb{I}_s} (|s|-t)\,dt + \int_{\widehat{\mathbb{I}}_s} e^{|s|-t}\,dt \\ &= -c_\varepsilon\left(|s||s|_\varepsilon - s^2 + \frac{s^2}{2} - \frac{(|s|_\varepsilon)^2}{2}\right) + e^{|s|}\int_{\widehat{\mathbb{I}}_s} e^{-t}\,dt \\ &= \frac{e^{-\varepsilon}\varepsilon}{2} + e^{-\varepsilon} < 2 = \alpha p(s),\end{aligned}$$

$$\begin{aligned}\int \boldsymbol{K}(s,t)p(s)\,ds &= -c_\varepsilon\left(\int_{\mathbb{I}_t^+} + \int_{\mathbb{I}_t^-}\right)(|s|-t)\,ds + \chi_{[\varepsilon,\infty)}(t)\int_{\mathbb{I}_t^\pm} e^{|s|-t}\,ds \\ &= -c_\varepsilon(t^2 - t_\varepsilon^2 - 2t^2 + 2tt_\varepsilon) + \chi_{[\varepsilon,\infty)}(t)e^{-t}\int_{\mathbb{I}_t^\pm} e^{|s|}\,ds \\ &= c_\varepsilon(t-t_\varepsilon)^2 + 2\chi_{[\varepsilon,\infty)}(t)e^{-t}e^{t-\varepsilon} \leqslant e^{-\varepsilon}\varepsilon + 2e^{-\varepsilon} < 2 = \beta q(t);\end{aligned}$$

here $p(s) = q(t) = 1$ for all $s$ and $t$ in $\mathbb{R}$ and $\alpha = \beta = 2$.

The function $\boldsymbol{K}$ in (2.20) is obviously continuous at each point $(s,t)$, but it fails to belong to $C(\mathbb{R}^2, \mathbb{C})$, because if $\alpha < 0$ is fixed, then $\boldsymbol{K}(s,t)$ takes the same non-zero value ($-c_\varepsilon\alpha$ or $e^\alpha$) for all points $(s,t)$ in $\mathbb{R}^2$ for which $|s| - t = \alpha$.

Write $f(t) := t$, $g(t) := e^{-t}$, $f^*(s) := |s|$, and $g^*(s) := e^{|s|}$, and let $\chi_A$ denote the characteristic function of the set $A$. With these notations the two Carleman functions $\boldsymbol{k}$, $\boldsymbol{k}'\colon \mathbb{R} \to L^2$ associated with $\boldsymbol{K}$ can then be written as

$$\begin{aligned}\boldsymbol{k}(s) = \overline{\boldsymbol{K}(s,\cdot)} &= -c_\varepsilon(|s| - f(\cdot))\chi_{\mathbb{I}_s}(\cdot) + e^{|s|}\chi_{\widehat{\mathbb{I}}_s}(\cdot)g(\cdot), \\ \boldsymbol{k}'(t) = \boldsymbol{K}(\cdot,t) &= \chi_{[0,\infty)}(t)(-c_\varepsilon(f^*(\cdot) - t)\chi_{\mathbb{I}_t^+\cup\mathbb{I}_t^-}(\cdot) + \chi_{[\varepsilon,\infty)}(t)e^{-t}\chi_{\mathbb{I}_t^\pm}(\cdot)g^*(\cdot)),\end{aligned}$$

with $s$ and $t$ running over $\mathbb{R}$. What is being asserted is that the functions $\boldsymbol{k}$ and $\boldsymbol{k}'$ are continuous from $\mathbb{R}$ into $L^2$, but they both are not in $C(\mathbb{R}, L^2)$, failing to satisfy (2.1). Indeed, if $s \to s_0 \in \mathbb{R}$ and $0 \leqslant t \to t_0 \in \mathbb{R}$, then

$$\begin{aligned}\|\boldsymbol{k}(s) - \boldsymbol{k}(s_0)\| &\leqslant c_\varepsilon\left(\left\||s|\chi_{\mathbb{I}_s} - |s_0|\chi_{\mathbb{I}_{s_0}}\right\| + \left\|f\left(\chi_{\mathbb{I}_{s_0}} - \chi_{\mathbb{I}_s}\right)\right\|\right) \\ &\quad + \left\|\left(e^{|s|}\chi_{\widehat{\mathbb{I}}_s} - e^{|s_0|}\chi_{\widehat{\mathbb{I}}_{s_0}}\right)g\right\| \\ &\leqslant c_\varepsilon\left(||s| - |s_0||\,\|\chi_{\mathbb{I}_s}\| + |s_0|\left\|\chi_{\mathbb{I}_s\Delta\mathbb{I}_{s_0}}\right\| + \left\|f\chi_{\mathbb{I}_s\Delta\mathbb{I}_{s_0}}\right\|\right) \\ &\quad + |e^{|s|} - e^{|s_0|}|\left\|\chi_{\widehat{\mathbb{I}}_s}g\right\| + e^{|s_0|}\left\|\chi_{\widehat{\mathbb{I}}_s\Delta\widehat{\mathbb{I}}_{s_0}}g\right\| \\ &\leqslant c_\varepsilon\left(\sqrt{\varepsilon}||s| - |s_0|| + (|s_0| + \max\{|s|,|s_0|\} + \varepsilon)\sqrt{2||s| - |s_0||}\right) \\ &\quad + \frac{1}{\sqrt{2}}|e^{|s|} - e^{|s_0|} + e^{|s_0|}\sqrt{2||s| - |s_0||} \to 0\end{aligned}$$

and

$$\begin{aligned}\|\boldsymbol{k'}(t)-\boldsymbol{k'}(t_0)\| &\leqslant c_\varepsilon\left(\left\|f^*\chi_{(\mathbb{I}_t^+\cup\mathbb{I}_t^-)\Delta(\mathbb{I}_{t_0}^+\cup\mathbb{I}_{t_0}^-)}\right\|+\left\|t\chi_{\mathbb{I}_t^+\cup\mathbb{I}_t^-}-t_0\chi_{\mathbb{I}_{t_0}^+\cup\mathbb{I}_{t_0}^-}\right\|\right)\\ &\quad+\left\|\left(\chi_{[\varepsilon,\infty)}(t)e^{-t}\chi_{\mathbb{I}_t^\pm}-\chi_{[\varepsilon,\infty)}(t_0)e^{-t_0}\chi_{\mathbb{I}_{t_0}^\pm}\right)g^*\right\|\\ &\leqslant c_\varepsilon\left(\|f^*\chi_{(\mathbb{I}_t^+\cup\mathbb{I}_t^-)\Delta(\mathbb{I}_{t_0}^+\cup\mathbb{I}_{t_0}^-)}\right\|+|t-t_0|\left\|\chi_{\mathbb{I}_t^+\cup\mathbb{I}_t^-}\|+t_0\left\|\chi_{(\mathbb{I}_t^+\cup\mathbb{I}_t^-)\Delta(\mathbb{I}_{t_0}^+\cup\mathbb{I}_{t_0}^-)}\right\|\right)\\ &\quad+|\chi_{[\varepsilon,\infty)}(t)e^{-t}-\chi_{[\varepsilon,\infty)}(t_0)e^{-t_0}|\left\|\chi_{\mathbb{I}_t^\pm}g^*\right\|+\chi_{[\varepsilon,\infty)}(t_0)e^{-t_0}\left\|\chi_{\mathbb{I}_t^\pm\Delta\mathbb{I}_{t_0}^\pm}g^*\right\|\to 0,\end{aligned}$$

respectively. But,

$$\begin{aligned}\|\boldsymbol{k}(s)\|^2 &= c_\varepsilon^2\,\|(|s|-f)\chi_{\mathbb{I}_s}\|^2+e^{2|s|}\left\|\chi_{\widehat{\mathbb{I}}_s}g\right\|^2\\ &= c_\varepsilon^2\int_{\mathbb{I}_s}(|s|-t)^2\,dt+e^{2|s|}\int_{\widehat{\mathbb{I}}_s}e^{-2t}\,dt\\ &= c_\varepsilon^2\frac{\varepsilon^3}{3}+\frac{e^{-2\varepsilon}}{2}\not\to 0\quad\text{as } |s|\to\infty,\end{aligned}$$

$$\begin{aligned}\|\boldsymbol{k'}(t)\|^2 &= c_\varepsilon^2\left\|(f^*-t)\chi_{\mathbb{I}_t^+\cup\mathbb{I}_t^-}\right\|^2+\chi_{[\varepsilon,\infty)}(t)e^{-2t}\left\|\chi_{\mathbb{I}_t^\pm}g^*\right\|^2\\ &= c_\varepsilon^2\int_{\mathbb{I}_t^+}(|s|-t)^2\,ds+c_\varepsilon^2\int_{\mathbb{I}_t^-}(|s|-t)^2\,ds+\chi_{[\varepsilon,\infty)}(t)e^{-2t}\int_{\mathbb{I}_t^\pm}e^{2|s|}\,ds\\ &= 2c_\varepsilon^2\frac{(t-t_\varepsilon)^3}{3}+\chi_{[\varepsilon,\infty)}(t)(e^{-2\varepsilon}-e^{-2t})\not\to 0\quad\text{as } t\to+\infty.\end{aligned}$$

What follows is a brief discussion of some direct consequences of Definition 1 relevant for this paper. In the first place, note that the conditions figuring in Definition 1 do not depend on each other in general; it is therefore natural to discuss the role played by each of them separately. The more restrictive of these conditions is (i), in the sense that it rules out the possibility for any $K^0$-kernel (unless that kernel is identically zero) of being a function depending only on the sum, difference, or product of the independent variables $s$, $t$; there are many other less trivial examples of inadmissible dependences (see also Example 1). These restrictions may be of use in constructing examples of bi-Carleman kernels which have both the properties (ii) and (iii), but do not enjoy (i); for another reason of existence of such a type of bi-Carleman kernels, we refer to the remark made in [36, p. 115] regarding boundedly supported kernels. Conversely it can be asserted, e.g., that if a function $\boldsymbol{T}\in C(\mathbb{R}^2,\mathbb{C})$ additionally satisfies for all $(s,t)\in\mathbb{R}^2$ the inequality $|\boldsymbol{T}(s,t)|\leqslant p(s)q(t)$, with $p$, $q$ being non-negative $C(\mathbb{R},\mathbb{R})$-functions square integrable over $\mathbb{R}$, then $\boldsymbol{T}$ is a $K^0$-kernel, i.e., the Carleman functions, $\boldsymbol{t}$ and $\boldsymbol{t'}$, it yields are both in $C(\mathbb{R},L^2)$. This assertion, obviously, pertains only to Hilbert-Schmidt kernels, and may be proved by an extension from the positive definite case with $p(s)\equiv q(s)\equiv(\boldsymbol{T}(s,s))^{\frac{1}{2}}$ to this general case of Buescu's argument in [2, pp. 247–249].

A few remarks are in order here concerning what can immediately be inferred from the $C(\mathbb{R},L^2)$-behaviour of the Carleman functions, $\boldsymbol{t}$ and $\boldsymbol{t'}$, associated with a given $K^0$-kernel, $\boldsymbol{T}$ (thought of as a kernel of an integral operator $T\in\mathfrak{R}(L^2)$):

1) The images of $\mathbb{R}$ under $\boldsymbol{t}$ and $\boldsymbol{t'}$, i.e.,

$$\boldsymbol{t}(\mathbb{R}):=\bigcup_{s\in\mathbb{R}}\boldsymbol{t}(s)\quad\text{and}\quad\boldsymbol{t'}(\mathbb{R}):=\bigcup_{s\in\mathbb{R}}\boldsymbol{t'}(s),\tag{2.21}$$

are relatively compact sets in $L^2$.

2) The *Carleman norm-functions* $\boldsymbol{\tau}$ and $\boldsymbol{\tau}'$, defined on $\mathbb{R}$ by $\boldsymbol{\tau}(s) = \|\boldsymbol{t}(s)\|$ and $\boldsymbol{\tau}'(s) = \|\boldsymbol{t}'(s)\|$, respectively, are continuous and vanish at infinity, i.e.,

$$\boldsymbol{\tau}, \boldsymbol{\tau}' \in C(\mathbb{R}, \mathbb{R}). \tag{2.22}$$

3) The images $Tf$ and $T^*f$ of any $f \in L^2$ under $T$ and $T^*$, respectively, have $C(\mathbb{R}, \mathbb{C})$-representatives in $L^2$, $[Tf]$ and $[T^*f]$, given pointwise by

$$[Tf](s) = \langle f, \boldsymbol{t}(s)\rangle, \quad [T^*f](s) = \langle f, \boldsymbol{t}'(s)\rangle \quad \text{at each } s \text{ in } \mathbb{R}. \tag{2.23}$$

4) Using (2.23), it is easy to deduce that $\boldsymbol{t}(\mathbb{R})^{\perp} = \operatorname{Ker} T$ and $\boldsymbol{t}'(\mathbb{R})^{\perp} = \operatorname{Ker} T^*$. (Indeed:

$$f \in \boldsymbol{t}(\mathbb{R})^{\perp} \iff \langle f, \boldsymbol{t}(s)\rangle = 0 \ \forall s \in \mathbb{R} \iff f \in \operatorname{Ker} T,$$
$$f \in \boldsymbol{t}'(\mathbb{R})^{\perp} \iff \langle f, \boldsymbol{t}'(s)\rangle \ \forall s \in \mathbb{R} \iff f \in \operatorname{Ker} T^*.)$$

The orthogonality between the range of an operator and the null-space of its adjoint then yields

$$\begin{aligned} \operatorname{Span}(\boldsymbol{t}(\mathbb{R})) &= \left(\boldsymbol{t}(\mathbb{R})^{\perp}\right)^{\perp} = \overline{\operatorname{Ran} T^*}, \\ \operatorname{Span}\left(\boldsymbol{t}'(\mathbb{R})\right) &= \left(\boldsymbol{t}'(\mathbb{R})^{\perp}\right)^{\perp} = \overline{\operatorname{Ran} T}. \end{aligned} \tag{2.24}$$

5) The $n$-th iterant $\boldsymbol{T}^{[n]}$ $(n \geqslant 2)$ of the $K^0$-kernel $\boldsymbol{T}$,

$$\boldsymbol{T}^{[n]}(s,t) := \underbrace{\int \cdots \int}_{n-1} \boldsymbol{T}(s, x_1) \ \ldots \ \boldsymbol{T}(x_{n-1}, t)\, dx_1 \ \ldots \ dx_{n-1} \left(= \langle T^{n-2}\left(\boldsymbol{t}'(t)\right), \boldsymbol{t}(s)\rangle\right), \tag{2.25}$$

is a $K^0$-kernel that defines the integral operator $T^n$. More generally, every two $K^0$-kernels $\boldsymbol{P}$ and $\boldsymbol{Q}$ might be said to be *multipliable* with each other, in the sense that their convolution

$$\boldsymbol{C}(s,t) := \int \boldsymbol{P}(s,x)\boldsymbol{Q}(x,t)\, dx \left(= \langle \boldsymbol{q}'(t), \boldsymbol{p}(s)\rangle\right)$$

exists at every point $(s,t) \in \mathbb{R}^2$, and forms a $K^0$-kernel that defines the product operator $C = PQ$. Indeed,

$$\begin{aligned} \int \langle \boldsymbol{q}'(t), \boldsymbol{p}(s)\rangle f(t)\, dt &= \int \left(\int \boldsymbol{P}(s,x)\boldsymbol{Q}(x,t)\, dx\right) f(t)\, dt = \langle f, Q^*(\boldsymbol{p}(s))\rangle \\ &= \langle Qf, \boldsymbol{p}(s)\rangle = \int \boldsymbol{P}(s,x) \left(\int \boldsymbol{Q}(x,t) f(t)\, dt\right) dx = [PQf](s) \end{aligned} \tag{2.26}$$

for every $f$ in $L^2$ and every $s$ in $\mathbb{R}$. Since both $\boldsymbol{p}$ and $\boldsymbol{q}'$ are in $C(\mathbb{R}, L^2)$ and both $P$ and $Q$ are in $\mathfrak{R}(L^2)$, the fact that $\boldsymbol{C}$ satisfies Definition 1 may be derived from the joint continuity of the inner product in its two arguments, which helps in proving (i), and from Proposition 1, whereby

$$\boldsymbol{c}(s) = \overline{\boldsymbol{C}(s,\cdot)} = Q^*(\boldsymbol{p}(s)) \quad \text{and} \quad \boldsymbol{c}'(s) = \boldsymbol{C}(\cdot, s) = P\left(\boldsymbol{q}'(s)\right) \quad \text{for every } s \text{ in } \mathbb{R},$$

which helps in proving both (ii) and (iii).

2.7. **Sub-$\boldsymbol{K^0}$-kernels.** Let $\boldsymbol{T}$ be a $K^0$-kernel defining an integral operator $T$, and impose on $\boldsymbol{T}$ an extra condition of being of ad hoc parquet-like support. Namely, let

(iv) there exist positive real numbers $\tau_n$ $(n \in \mathbb{N})$ strictly increasing with $n$ to infinity, $\tau_n \uparrow \infty$ as $n \to \infty$, such that for each fixed $n$, the *subkernels of* $\boldsymbol{T}$, $\boldsymbol{T}_n$ and $\widetilde{\boldsymbol{T}}_n$, defined on $\mathbb{R}^2$ by

$$\boldsymbol{T}_n(s,t) = \chi_n(s)\boldsymbol{T}(s,t) \quad \text{and} \quad \widetilde{\boldsymbol{T}}_n(s,t) = \boldsymbol{T}_n(s,t)\chi_n(t), \tag{2.27}$$

are $K^0$-kernels, and the integral operators,

$$T_n := P_n T \quad \text{and} \quad \widetilde{T}_n := P_n T P_n, \tag{2.28}$$

they induce on $L^2$ are nuclear.

Here in (iv), as in the rest of the paper, $\chi_n$ stands for the characteristic function of the open interval $\mathbb{I}_n = (-\tau_n, \tau_n)$, and $P_n$ for an orthogonal projection of $\mathfrak{P}(L^2)$ defined on each $f \in L^2$ by $P_n f = \chi_n f$; hence $(I - P_n)f = \widehat{\chi}_n f$ for each $f \in L^2$, with $\widehat{\chi}_n$ standing for the characteristic function of the set $\widehat{\mathbb{I}}_n = \mathbb{R} \setminus \mathbb{I}_n$.

A $K^0$-kernel $\boldsymbol{T}$ which satisfies condition (iv) necessarily has to vanish everywhere on the straight lines $s = \pm\tau_n$ and $t = \pm\tau_n$, parallel to the $t$ and $s$ axes, respectively; record this "parquet" property by writing

$$\boldsymbol{\tau}(\pm\tau_n) = \boldsymbol{\tau}'(\pm\tau_n) = \boldsymbol{T}(\pm\tau_n, t) = \boldsymbol{T}(s, \pm\tau_n) = 0 \quad \text{for all } s,\, t \text{ in } \mathbb{R} \text{ and } n \text{ in } \mathbb{N}. \tag{2.29}$$

The $P_n$ ($n \in \mathbb{N}$) clearly form a sequence of orthogonal projections increasing to $I$ with respect to the strong operator topology, so that, for every $f \in L^2$,

$$\|(P_n - I)\, f\| \searrow 0 \quad \text{as } n \to \infty. \tag{2.30}$$

So it follows immediately from (2.28) that, as $n \to \infty$,

$$\begin{aligned} \|(T_n - T)f\| \to 0, \quad \|(\widetilde{T}_n - T)f\| \to 0, \\ \|(T_n^* - T^*)f\| \to 0, \quad \|(\widetilde{T}_n^* - T^*)f\| \to 0. \end{aligned} \tag{2.31}$$

Among the subkernels defined in (2.27), the $\boldsymbol{T}_n$'s have more in common with the original kernel $\boldsymbol{T}$, inasmuch as $[T_n f](s) = \int \boldsymbol{T}(s,t) f(t)\, dt$ for every $f$ in $L^2$ and every $s$ in $\mathbb{I}_n$, while the tilded subkernels $\widetilde{\boldsymbol{T}}_n$ are more suitable to deal with $\boldsymbol{T}$ being Hermitian, i.e., satisfying $\boldsymbol{T}(s,t) = \overline{\boldsymbol{T}(t,s)}$ for all $s$, $t \in \mathbb{R}$, because then they all are also Hermitian.

Now we list some basic properties of the subkernels defined in (2.27), most of which are obvious from the definition:

$$|\boldsymbol{T}_n(s,t)| \leqslant |\boldsymbol{T}(s,t)|, \quad |\widetilde{\boldsymbol{T}}_n(s,t)| \leqslant |\boldsymbol{T}(s,t)|, \quad \text{for all } s,\, t \in \mathbb{R}, \tag{2.32}$$

$$\lim_{n\to\infty} \|\boldsymbol{T}_n - \boldsymbol{T}\|_{C(\mathbb{R}^2,\mathbb{C})} = 0, \quad \lim_{n\to\infty} \|\widetilde{\boldsymbol{T}}_n - \boldsymbol{T}\|_{C(\mathbb{R}^2,\mathbb{C})} = 0, \tag{2.33}$$

$$\int\!\!\int |\boldsymbol{T}_n(s,t)|^2\, dt\, ds < \infty, \quad \int\!\!\int |\widetilde{\boldsymbol{T}}_n(s,t)|^2\, dt\, ds < \infty, \tag{2.34}$$

$$\begin{aligned} \lim_{n\to\infty} \|\boldsymbol{t}_n - \boldsymbol{t}\|_{C(\mathbb{R},L^2)} = 0, \quad \lim_{n\to\infty} \|\boldsymbol{t}'_n - \boldsymbol{t}'\|_{C(\mathbb{R},L^2)} = 0, \\ \lim_{n\to\infty} \|\widetilde{\boldsymbol{t}}_n - \boldsymbol{t}\|_{C(\mathbb{R},L^2)} = 0, \quad \lim_{n\to\infty} \|\widetilde{\boldsymbol{t}}'_n - \boldsymbol{t}'\|_{C(\mathbb{R},L^2)} = 0, \end{aligned} \tag{2.35}$$

where, for each $n$ in $\mathbb{N}$, the Carleman functions $\boldsymbol{t}_n$, $\boldsymbol{t}'_n$ and $\widetilde{\boldsymbol{t}}_n$, $\widetilde{\boldsymbol{t}}'_n$, associated to the subkernels $\boldsymbol{T}_n$ and $\widetilde{\boldsymbol{T}}_n$, are defined, as usual, to be

$$\begin{aligned} \boldsymbol{t}_n(s) = \overline{\boldsymbol{T}_n(s,\cdot)} \ (= \chi_n(s)\boldsymbol{t}(s)), \quad \boldsymbol{t}'_n(t) = \boldsymbol{T}_n(\cdot,t) \ (= P_n\left(\boldsymbol{t}'(t)\right)), \\ \widetilde{\boldsymbol{t}}_n(s) = \overline{\widetilde{\boldsymbol{T}}_n(s,\cdot)} \ (= \chi_n(s) P_n\left(\boldsymbol{t}(s)\right)), \quad \widetilde{\boldsymbol{t}}'_n(t) = \widetilde{\boldsymbol{T}}_n(\cdot,t) \ (= \chi_n(t) P_n\left(\boldsymbol{t}'(t)\right)) \end{aligned} \tag{2.36}$$

for every $s$, $t \in \mathbb{R}$. The limits in (2.35) all hold due to (ii), (iii), (2.30), and a result from [14, Lemma 3.7, p. 151] (or [5, Theorem 3.2, p. 124]). The result, just referred to, will be used in the text so often that it should be explicitly stated:

**Lemma 1.** *Let $B_n$, $B \in \mathfrak{R}(L^2)$, and suppose that for every $f \in L^2$, $\|B_n f - Bf\| \to 0$ as $n \to \infty$. Then for any relatively compact set $U$ in $L^2$,*

$$\sup_{f\in U} \|B_n f - Bf\| \to 0 \quad \text{as } n \to \infty. \tag{2.37}$$

Applying this lemma to the sets $\boldsymbol{t}(\mathbb{R})$ and $\boldsymbol{t}'(\mathbb{R})$ in (2.21) immediately gives that, as $n \to \infty$,

$$\sup_{s\in\mathbb{R}} \|(B_n - B)(\boldsymbol{t}(s))\| \to 0 \quad \text{and} \quad \sup_{t\in\mathbb{R}} \|(B_n - B)(\boldsymbol{t}'(t))\| \to 0, \tag{2.38}$$

and putting $B_n = P_n$, $B = I$ here causes the three last limits in (2.35) to be zero.

If, in the lemma which follows, the operators $B_n$ are each equal to $I$, then its conclusion for that case can be easily derived from (2.35) and Theorem 1. However, for the general case, we supply an independent proof based only on Definition 1 and condition (iv).

**Lemma 2.** *For each $n$ in $\mathbb{N}$, let $B_n \in \mathfrak{R}(L^2)$ and define two functions $\boldsymbol{b}_n$, $\boldsymbol{b}'_n \colon \mathbb{R} \to L^2$ by $\boldsymbol{b}_n(s) = B_n(\boldsymbol{t}_n(s))$, $\boldsymbol{b}'_n(t) = B_n^*(\boldsymbol{t}'_n(t))$. Then, if $\sup_{n\in\mathbb{N}} \|B_n\| < \infty$, each of the subsets*

$$\mathfrak{b} = \bigcup_{n=1}^{\infty} \boldsymbol{b}_n \subset C(\mathbb{R}, L^2) \quad \text{and} \quad \mathfrak{b}' = \bigcup_{n=1}^{\infty} \boldsymbol{b}'_n \subset C(\mathbb{R}, L^2) \tag{2.39}$$

*is bounded in $C(\mathbb{R}, L^2)$, and is equicontinuous on $\mathbb{R} \cup \{\infty\}$, i.e., has property* (A2) *of Theorem 1, for $k = 1$ and $B = L^2$.*

*Proof.* Let $M > 0$ be that constant for which $\|B_n\| \leqslant M$ for all $n$ in $\mathbb{N}$. Fix any $n$ in $\mathbb{N}$, and make use of property (2.29) and formulae (2.36) to establish the following inequalities, valid for all $x$ and $y$ in $\mathbb{R}$:

$$\begin{gathered}\|\boldsymbol{b}_n(x) - \boldsymbol{b}_n(y)\| \leqslant \|B_n\| \|\boldsymbol{t}_n(x) - \boldsymbol{t}_n(y)\| \\ \leqslant \begin{cases} M \cdot \max\{\|\boldsymbol{t}(x)\|, \|\boldsymbol{t}(y)\|\} & \text{if } |x-y| = |\tau_n - x| + |\tau_n - y|, \\ M \cdot \max\{\|\boldsymbol{t}(x)\|, \|\boldsymbol{t}(y)\|\} & \text{if } |x-y| = |\tau_n + x| + |\tau_n + y|, \\ M\|\boldsymbol{t}(x) - \boldsymbol{t}(y)\| & \text{otherwise}, \end{cases} \\ \|\boldsymbol{b}'_n(x) - \boldsymbol{b}'_n(y)\| \leqslant \|B_n^* P_n\| \|\boldsymbol{t}'(x) - \boldsymbol{t}'(y)\| \leqslant M\|\boldsymbol{t}'(x) - \boldsymbol{t}'(y)\|.\end{gathered} \tag{2.40}$$

In particular, when written for a fixed $y$ satisfying $\boldsymbol{\tau}(y) = \boldsymbol{\tau}'(y) = 0$, these inequalities lead to the following estimates, valid for all $x$ in $\mathbb{R}$ and $n$ in $\mathbb{N}$:

$$\|\boldsymbol{b}_n(x)\| \leqslant M\|\boldsymbol{\tau}\|_{C(\mathbb{R},\mathbb{R})}, \quad \|\boldsymbol{b}'_n(x)\| \leqslant M\|\boldsymbol{\tau}'\|_{C(\mathbb{R},\mathbb{R})}, \tag{2.41}$$

thereby proving the boundedness of $\mathfrak{b}$ and $\mathfrak{b}'$ in $C(\mathbb{R}, L^2)$. Since, by Definition 1, both $\boldsymbol{t}$ and $\boldsymbol{t}'$ are uniformly continuous vanishing at infinity, it follows that given any $\varepsilon > 0$, a $\delta > 0$ can be found such that for $x$ and $y$ in $\mathbb{R}$,

$$\max\{\|\boldsymbol{t}(x) - \boldsymbol{t}(y)\|, \|\boldsymbol{t}'(x) - \boldsymbol{t}'(y)\|\} < \varepsilon$$

whenever $|x - y| < \delta$ or/and whenever both $|x| > \dfrac{1}{\delta}$ and $\boldsymbol{\tau}(y) = \boldsymbol{\tau}'(y) = 0$. Combining this with (2.40) gives that for all the $x$ and $y$ values so selected,

$$\max\{\|\boldsymbol{b}_n(x) - \boldsymbol{b}_n(y)\|, \|\boldsymbol{b}'_n(x) - \boldsymbol{b}'_n(y)\|\} < M\varepsilon \quad \text{for every } n \in \mathbb{N}. \tag{2.42}$$

Hence, each of the subsets $\mathfrak{b}$ and $\mathfrak{b}'$ in (2.39) is equicontinuous on $\mathbb{R} \cup \{\infty\}$, in the sense of condition (A2) of Theorem 1. The lemma is proved. □

*Remark* 3. We should mention the papers [34], [6], and [35]), where conditions analogous to (2.32)-(2.34) and (2.40) occur. In more detail, under consideration in these works is the bi-Carleman kernel $\boldsymbol{K}$, supported in $(0,1)^2$, which is merely postulated to be approximable, in the pointwise sense that for almost every $(s,t) \in (0,1)^2$

$$\boldsymbol{K}(s,t) = \lim_{m\to\infty} \boldsymbol{K}_m(s,t) \text{ (compare with (2.33))}, \tag{2.43}$$

by the kernels $\boldsymbol{K}_m$ $(m \in \mathbb{N})$, each of which possesses all the following properties:

$$|\boldsymbol{K}_m(s,t)| \leqslant |\boldsymbol{K}(s,t)| \text{ (compare with (2.32))}, \tag{2.44}$$

$$\int_0^1 \int_0^1 |\boldsymbol{K}_m(s,t)|^2 \, ds \, dt < \infty \text{ (compare with (2.34))}, \tag{2.45}$$

$$\int_0^1 |\boldsymbol{K}_m(s,t) - \boldsymbol{K}_m(s_0,t)|^2 \, dt \leqslant (\sigma(s,s_0))^2 \text{ (compare with (2.40))} \tag{2.46}$$

$$\int_0^1 |\boldsymbol{K}_m^*(s,t) - \boldsymbol{K}_m^*(s_0,t)|^2 \, dt \leqslant (\sigma(s,s_0))^2 \text{ (compare with (2.40))}, \tag{2.47}$$

$$\int_0^1 \boldsymbol{K}_m(s,x)\boldsymbol{K}_m^*(x,t) \, dx = \int_0^1 \boldsymbol{K}_m^*(s,x)\boldsymbol{K}_m(x,t) \, dx, \tag{2.48}$$

where the function $\sigma\colon (0,1)^2 \to \mathbb{R}$ is independent of $m$ and such that $\sigma(s,s_0)$ tends to zero as $s_0$ tends to $s$ for almost every $s \in (0,1)$, and where $\boldsymbol{K}_m^*(s,t) = \overline{\boldsymbol{K}_m(t,s)}$. The listed properties of the approximating kernels $\boldsymbol{K}_m$ are, in fact, an abstraction of those of the subkernels used by Carleman in his famous work [4]. There he investigates real-valued kernel $\boldsymbol{K}$ defined on $(0,1)^2$ and satisfying the following three conditions:

(C1) $\boldsymbol{K}(s,t) = \boldsymbol{K}(t,s)$ for all $s$, $t \in (0,1)$;

(C2) $\boldsymbol{k}(s)^2 = \int_0^1 \boldsymbol{K}(s,t)^2 \, dt < \infty$ for all $s \in (0,1)$, and the limit

$$\lim_{s \to s_0} \int_0^1 (\boldsymbol{K}(s_0,t) - \boldsymbol{K}(s,t))^2 \, dt = 0$$

holds for all $s_0$ in $(0,1)$ except possibly for an at most countable subset $\{s_i\} \subset (0,1)$ having at most a finite number of accumulation points;

(C3) In the set $\{s_i\}$, there is a finite number of points $\eta_1, \eta_2, \ldots, \eta_l$ such that for all sufficiently small $\varepsilon > 0$,

$$\int_{\mathbb{I}_\varepsilon} \boldsymbol{k}(s)^2 \, ds < \infty \quad \text{with } \mathbb{I}_\varepsilon := (0,1) \setminus \left( \bigcup_{k=1}^{l} \{s \in (0,1) \mid |s - \eta_k| < \varepsilon\} \right).$$

Fix a sequence $\{\varepsilon_m\}_{m=1}^\infty$ of sufficiently small positive numbers decreasing to zero. For each $m$ in $\mathbb{N}$, let $\chi_m$ be the characteristic function of the set $\mathbb{I}_{\varepsilon_m}$, and define the subkernel $\boldsymbol{K}_m(s,t) := \chi_m(s)\boldsymbol{K}(s,t)\chi_m(t)$ (compare (2.27)). It is not hard to check that the above conditions (2.43)-(2.48) are satisfied by all such definable subkernels $\boldsymbol{K}_m$ $(m \in \mathbb{N})$.

Every $K^0$-kernel $\boldsymbol{T}$, obviously, satisfies Carleman's conditions (C2)-(C3) in $s$ (and in $t$) with only one exceptional point $s_1 = \eta_1 = \infty$, and hence with $\mathbb{I}_\varepsilon$ replaced by the set $\{s \mid |s| < \frac{1}{\varepsilon}\}$. We shall deal with just a small part of Carleman's [4] theory in Subsection 4.4 of the paper, where we shall try to construct spectral functions for Hermitian $K^0$-kernels as limits of sequences of spectral functions for their subkernels.

2.8. **Unitary reductions to $\boldsymbol{K^0}$-kernels.** We would like to close this Section 2 with a unitary equivalence result which is essentially contained in Theorem of [24], where it is proved for $K^0$-kernels supported in the quarter plane of the first quadrant (i.e., for operators acting from $L^2[0,\infty)$ into itself) and with the explicit description of the real sequence $\{\tau_n\}$ (denoted there by "$\{t_n\}$") needed to enforce the condition (iv). (See also Theorem 1 in [23] and in [25].) Here is a whole-line variant of that result, which will be the basic ingredient for our approach to deal with equations like (1.1) (see also Remark 4 below):

**Proposition 2.** *Suppose that $S$ is a bi-integral operator on $L^2$. Then there exists a unitary operator $U\colon L^2 \to L^2$ such that the operator $T = USU^{-1}$ is a bi-Carleman operator on $L^2$, whose kernel is a $K^0$-kernel satisfying condition* (iv) *of Subsection 2.7.*

The (constructive) proof for this proposition can, with a slight modification, be read off the proof of Theorem 9 in [27, Section 4]. The modification that has to be made is to use, instead of the Lemarié-Meyer wavelet basis $\{u_n\}$ utilized in that proof, an orthonormal basis of $L^2$ formed by the compactly supported continuous functions of $s$:

$$u_{n_{r,l}}(s) = \begin{cases} \dfrac{\sqrt{2}}{l^2}\chi_{[\tau_{l-1},\tau_l]}(s)\sin\left(\dfrac{r\pi(s-\tau_{l-1})}{\tau_l-\tau_{l-1}}\right) & \text{if } l>0, \\ \dfrac{\sqrt{2}}{|l-1|^2}\chi_{[-\tau_{|l|+1},-\tau_{|l|}]}(s)\sin\left(\dfrac{r\pi(s+\tau_{|l|+1})}{\tau_{|l|+1}-\tau_{|l|}}\right) & \text{if } l\le 0, \end{cases} \tag{2.49}$$

where $r$ ranges through $\mathbb{N}$ and $l$ through the integers $\mathbb{Z}$; $\mathbb{N} = \cup_{l\in\mathbb{Z}}\{n_{r,l}\}_{r=1}^{\infty}$ is a disjoint decomposition of $\mathbb{N}$ into countably many sequences such that $n_{r,l} < n_{r+1,l}$ for all $r$ and all $l$; $\chi_{[a,b]}$ denotes the characteristic function of the interval $[a,b]$; $\tau_0 = 0$ and for each $n$ in $\mathbb{N}$, $\tau_n = 1^4 + \cdots + n^4$. (It is intended that these $\tau_n$'s are to have the role of those involved in condition (iv).) Accordingly, the allowed range of the values taken by the integers $i$ and $j$, destined to play the role of orders of differentiation in that proof, has to be restricted form $\mathbb{N}\cup\{0\}$ to $\{0\}$.

*Remark* 4. By virtue of Proposition 2 and Remark 2, we may and will confine our investigation with no loss of generality to second-kind integral equations (1.1) in which the kernel $\boldsymbol{T}$, the *original kernel* as we shall continue to call it below, induces a bi-Carleman operator $T \in \mathfrak{R}(L^2)$ and possesses all the properties (i)-(iv) (defined in Definition 1 and Subsection 2.7). These assumptions on $\boldsymbol{T}$ will remain in force for the rest of the paper, though not all of them will be employed in specific situations; the notations given in condition (iv) will therefore be used frequently without warning in the subsequent sections.

## 3. Resolvent $K^0$-kernels

3.1. **Resolvent kernels for $\boldsymbol{K^0}$-kernels.** We start this section with the definition of the resolvent kernel for a $K^0$-kernel, which differs from the general one given at the beginning of the introduction at least in that it does not make explicit use of the concept of Fredholm resolvent.

**Definition 2.** Let $\boldsymbol{T}$ be a $K^0$-kernel, let $\lambda$ be a complex number, and suppose that a $K^0$-kernel, to be denoted by $\boldsymbol{T}_{|\lambda}$, satisfies, for all $s$ and $t$ in $\mathbb{R}$, the two simultaneous integral equations

$$\boldsymbol{T}_{|\lambda}(s,t) - \lambda\int \boldsymbol{T}(s,x)\boldsymbol{T}_{|\lambda}(x,t)\,dx = \boldsymbol{T}(s,t), \tag{3.1}$$

$$\boldsymbol{T}_{|\lambda}(s,t) - \lambda\int \boldsymbol{T}_{|\lambda}(s,x)\boldsymbol{T}(x,t)\,dx = \boldsymbol{T}(s,t), \tag{3.2}$$

and the condition that for any $f$ in $L^2$,

$$\int\left|\int \boldsymbol{T}_{|\lambda}(s,t)f(t)\,dt\right|^2 ds < \infty. \tag{3.3}$$

Then the $K^0$-kernel $\boldsymbol{T}_{|\lambda}$ will be called the *resolvent kernel* for $\boldsymbol{T}$ at $\lambda$, and the functions $\boldsymbol{t}_{|\lambda}$ and $\boldsymbol{t}'_{|\lambda}$ of $C(\mathbb{R},L^2)$, defined via $\boldsymbol{T}_{|\lambda}$ by $\boldsymbol{t}_{|\lambda}(s) = \overline{\boldsymbol{T}_{|\lambda}(s,\cdot)}$ and $\boldsymbol{t}'_{|\lambda}(t) = \boldsymbol{T}_{|\lambda}(\cdot,t)$, will be called the *resolvent Carleman functions* for $\boldsymbol{T}$ at $\lambda$.

**Theorem 2.** *Let $T \in \mathfrak{R}(L^2)$ be an integral operator, with a kernel $\boldsymbol{T}$ that is a $K^0$-kernel, and let $\lambda$ be a complex number. Then* (a) *if $\lambda$ is a regular value for $T$, then the resolvent kernel for $\boldsymbol{T}$ exists at $\lambda$, and is a kernel of the Fredholm resolvent of $T$ at $\lambda$, i.e., $\left(T_{|\lambda}f\right)(s) = \int \boldsymbol{T}_{|\lambda}(s,t)f(t)\,dt$ for every $f$ in $L^2$ and almost every $s$ in $\mathbb{R}$;* (b) *if the resolvent kernel for $\boldsymbol{T}$ exists at $\lambda$, then $\lambda$ is a regular value for $T$.*

*Proof.* To prove statement (a), let $\lambda$ be an arbitrary but fixed regular value for $T$ ($\lambda \in \Pi(T)$), and define two functions $\boldsymbol{a}$, $\boldsymbol{a'}\colon \mathbb{R} \to L^2$ by writing

$$\boldsymbol{a}(s) = \left(\bar{\lambda}(T_{|\lambda})^* + I\right)(\boldsymbol{t}(s)), \quad \boldsymbol{a'}(s) = \left(\lambda T_{|\lambda} + I\right)\left(\boldsymbol{t'}(s)\right) \tag{3.4}$$

whenever $s \in \mathbb{R}$. So defined, $\boldsymbol{a}$ and $\boldsymbol{a'}$ then belong to the space $C(\mathbb{R}, L^2)$, as $\boldsymbol{t}$ and $\boldsymbol{t'}$ (the Carleman functions associated to the $K^0$-kernel $\boldsymbol{T}$) are in $C(\mathbb{R}, L^2)$, and $T_{|\lambda}$ (the Fredholm resolvent of $T$ at $\lambda$) is in $\mathfrak{R}(L^2)$.

The functions $\boldsymbol{A}$, $\boldsymbol{A'}\colon \mathbb{R}^2 \to \mathbb{C}$, given by the formulae

$$\begin{aligned} \boldsymbol{A}(s,t) &= \lambda \left\langle \boldsymbol{t'}(t), \boldsymbol{a}(s)\right\rangle + \boldsymbol{T}(s,t), \\ \boldsymbol{A'}(s,t) &= \bar{\lambda}\overline{\langle \boldsymbol{a'}(s), \boldsymbol{t}(t)\rangle} + \overline{\boldsymbol{T}(t,s)}, \end{aligned} \tag{3.5}$$

then belong to the space $C(\mathbb{R}^2, \mathbb{C})$, due to the continuity of the inner product as a function from $L^2 \times L^2$ to $\mathbb{C}$. By using (3.4) it is also seen from (3.5) that these functions are conjugate transposes of each other, i.e., $\boldsymbol{A'}(s,t) = \overline{\boldsymbol{A}(t,s)}$ for all $s$, $t \in \mathbb{R}$. Simple manipulations involving formulae (3.5), (2.23), and (3.4) give rise to the following two strings of equations, satisfied at all points $s$ in $\mathbb{R}$ by each function $f$ in $L^2$:

$$\begin{aligned} \int \boldsymbol{A}(s,t)f(t)\,dt &= \lambda \int \left\langle \boldsymbol{t'}(t), \boldsymbol{a}(s)\right\rangle f(t)\,dt + \int \boldsymbol{T}(s,t)f(t)\,dt \\ &= \left\langle f, \bar{\lambda}T^*(\boldsymbol{a}(s)) + \boldsymbol{t}(s)\right\rangle = \langle f, \boldsymbol{a}(s)\rangle, \\ \int \boldsymbol{A'}(s,t)f(t)\,dt &= \bar{\lambda} \int \overline{\langle \boldsymbol{a'}(s), \boldsymbol{t}(t)\rangle} f(t)\,dt + \int \overline{\boldsymbol{T}(t,s)} f(t)\,dt \\ &= \left\langle f, \lambda T\left(\boldsymbol{a'}(s)\right) + \boldsymbol{t'}(s)\right\rangle = \left\langle f, \boldsymbol{a'}(s)\right\rangle. \end{aligned}$$

The equality of the extremes of each of these strings implies that $\overline{\boldsymbol{A}(s,\cdot)} \in \boldsymbol{a}(s)$ and $\boldsymbol{A}(\cdot,s) \in \boldsymbol{a'}(s)$ for every fixed $s$ in $\mathbb{R}$. Furthermore, the following relations hold whenever $f$ is in $L^2$:

$$\begin{aligned} \int \boldsymbol{A}(\cdot,t)f(t)\,dt &= \langle f, \boldsymbol{a}(\cdot)\rangle = \left\langle \left(\lambda T_{|\lambda} + I\right) f, \boldsymbol{t}(\cdot)\right\rangle \\ &= \langle R_\lambda(T)f, \boldsymbol{t}(\cdot)\rangle = (TR_\lambda(T)f)(\cdot) = \left(T_{|\lambda}f\right)(\cdot) \in L^2, \end{aligned} \tag{3.6}$$

showing that the Fredholm resolvent $T_{|\lambda}$ of $T$ at $\lambda$ is an integral operator on $L^2$, with the function $\boldsymbol{A}$ as its kernel (compare this with (3.3)).

The inner product when written in the integral form and the above observations about $\boldsymbol{A}$ allow the defining relationships for $\boldsymbol{A}$ and $\boldsymbol{A'}$ (see (3.5)) to be respectively written as the integral equations

$$\begin{aligned} \boldsymbol{A}(s,t) &= \lambda \int \boldsymbol{A}(s,x)\boldsymbol{T}(x,t)\,dx + \boldsymbol{T}(s,t), \\ \boldsymbol{A}(s,t) &= \lambda \int \boldsymbol{T}(s,x)\boldsymbol{A}(x,t)\,dx + \boldsymbol{T}(s,t), \end{aligned}$$

holding for all $s$, $t \in \mathbb{R}$. Together with (3.6), these imply that the $K^0$-kernel $\boldsymbol{A}$ is a resolvent kernel for $\boldsymbol{T}$ at $\lambda$ (in the sense of Definition 2).

To prove statement (b), let there exist a $K^0$-kernel $\boldsymbol{T}_{|\lambda}$ satisfying (3.1) through (3.3). It is to be proved that $\lambda$ belongs to $\Pi(T)$, i.e., that the operator $I - \lambda T$ is invertible. To this effect, therefore, remark first that the integral operator $A$ given by $(Af)(s) = \int \boldsymbol{T}_{|\lambda}(s,t)f(t)\,dt$ is bounded from $L^2$ into $L^2$, owing to condition (3.3) and to the Banach Theorem (see [10, p. 14]). Then, due to the multipliability property of $K^0$-kernels (see (2.26)), the kernel-function equations (3.1) and (3.2) give rise to the operator equalities $(I-\lambda T)A = T$ and $A(I-\lambda T) = T$, respectively. The latter are easily seen to be equivalent respectively to the following ones $(I-\lambda T)(I+\lambda A) = I$ and $(I+\lambda A)(I-\lambda T) = I$, which together imply that the operator $I-\lambda T$ is invertible with inverse $I+\lambda A$. The theorem is proved. □

*Remark* 5. The proof just given establishes that resolvent kernels in the sense of Definition 2 are in one-to-one correspondence with Fredholm resolvents. In view of this correspondence: (1) $\Pi(T)$ might as well be defined as the set of all those $\lambda \in \mathbb{C}$ at which the resolvent kernel in the sense of Definition 2 exists (thus, whenever $\boldsymbol{T}_{|\lambda}$, $\boldsymbol{t}_{|\lambda}$, or $\boldsymbol{t}'_{|\lambda}$ appear in what follows, it may and will always be understood that $\lambda$ belongs to $\Pi(T)$); (2) the resolvent kernel $\boldsymbol{T}_{|\lambda}$ for the $K^0$-kernel $\boldsymbol{T}$ at $\lambda$ might as well be defined as that $K^0$-kernel which induces $T_{|\lambda}$, the Fredholm resolvent at $\lambda$ of that integral operator $T$ whose kernel is $\boldsymbol{T}$. Using (2.3) and (2.19), the values of the resolvent Carleman functions for $\boldsymbol{T}$ at each fixed regular value $\lambda \in \Pi(T)$ can therefore be ascertained by writing

$$\boldsymbol{t}_{|\lambda}(\cdot) = R^*_{\lambda}(T)(\boldsymbol{t}(\cdot)), \quad \boldsymbol{t}'_{|\lambda}(\cdot) = R_{\lambda}(T)\left(\boldsymbol{t}'(\cdot)\right), \tag{3.7}$$

where $\boldsymbol{t}$ and $\boldsymbol{t}'$ are Carleman functions corresponding to $\boldsymbol{T}$ (compare with (3.4) via (2.4)). The resolvent kernel $\boldsymbol{T}_{|\lambda}$ for $\boldsymbol{T}$, in its turn, can be exactly recovered from the knowledge of the resolvent Carleman functions $\boldsymbol{t}_{|\lambda}$ and $\boldsymbol{t}'_{|\lambda}$ by the formulae

$$\overline{\boldsymbol{T}_{|\lambda}(s,t)} = \bar{\lambda}\left\langle \boldsymbol{t}_{|\lambda}(s), \boldsymbol{t}'(t)\right\rangle + \overline{\boldsymbol{T}(s,t)}, \tag{3.8}$$

$$\boldsymbol{T}_{|\lambda}(s,t) = \lambda\left\langle \boldsymbol{t}'_{|\lambda}(t), \boldsymbol{t}(s)\right\rangle + \boldsymbol{T}(s,t), \tag{3.9}$$

respectively (compare with (3.5)). Formulae (3.7)-(3.9) will be of constant use in the subsequent analysis.

*Remark* 6. When, as in (2.8), $\lambda$ is in $d(T) := \{\lambda \in \mathbb{C} \mid |\lambda| < r(T)\}$, the values of the resolvent kernel $\boldsymbol{T}_{|\lambda}$ can be analytically computed using nothing but the iterants of $\boldsymbol{T}$; the reasoning goes as follows:

$$\begin{aligned}
\boldsymbol{T}_{|\lambda}(s,t) &= \boldsymbol{T}(s,t) + \lambda\langle R_{\lambda}(T)(\boldsymbol{t}'(t)), \boldsymbol{t}(s)\rangle && \text{by (3.9) and (3.7)}\\
&= \boldsymbol{T}(s,t) + \lambda\langle\left(\sum_{n=0}^{\infty}\lambda^n T^n\right)(\boldsymbol{t}'(t)), \boldsymbol{t}(s)\rangle && \text{by (2.8)}\\
&= \boldsymbol{T}(s,t) + \sum_{n=0}^{\infty}\left\langle \lambda^{n+1}T^n(\boldsymbol{t}'(t)), \boldsymbol{t}(s)\right\rangle && \text{by (2.8)}\\
&= \sum_{n=1}^{\infty}\lambda^{n-1}\boldsymbol{T}^{[n]}(s,t), && \text{by (2.25)}
\end{aligned}$$

where $\boldsymbol{T}^{[1]} = \boldsymbol{T}$. The series in the last line may be called the Neumann expansion for $\boldsymbol{T}_{|\lambda}$. The estimate

$$\|\boldsymbol{T}^{[n]}\|_{C(\mathbb{R}^2,\mathbb{C})} = \sup_{(s,t)\in\mathbb{R}^2}\left|\langle T^{n-2}\boldsymbol{t}'(t), \boldsymbol{t}(s)\rangle\right| \leqslant \|\boldsymbol{\tau}'\|_{C(\mathbb{R},\mathbb{R})}\|\boldsymbol{\tau}\|_{C(\mathbb{R},\mathbb{R})}\|T^{n-2}\|$$

implies that the $C\left(\mathbb{R}^2,\mathbb{C}\right)$-convergence of this series to $\boldsymbol{T}_{|\lambda}$ is also uniform in $\lambda$, i.e., that the following limit holds true whenever $\mathfrak{K}$ is a compact subset of $d(T)$:

$$\lim_{m\to\infty}\sup_{\lambda\in\mathfrak{K}}\|\boldsymbol{T}_{|\lambda}-\sum_{n=1}^{m}\lambda^{n-1}\boldsymbol{T}^{[n]}\|_{C(\mathbb{R}^2,\mathbb{C})}=0. \tag{3.10}$$

This formula efficiently solves bi-Carleman version of Korotkov's problem (∗) (see the introduction) at $\lambda$ running throughout the disk $d(T)$; however, the latter is in general a proper subset of $\Pi(T)$.

3.2. **Resolvent kernels for sub-$\boldsymbol{K}^0$-kernels.** Here, as subsequently, we shall denote the resolvent kernels at $\lambda$ for the subkernel $\boldsymbol{T}_n$ (resp., $\widetilde{\boldsymbol{T}}_n$) by $\boldsymbol{T}_{n|\lambda}$ (resp., $\widetilde{\boldsymbol{T}}_{n|\lambda}$), and the resolvent Carleman functions for these subkernels at $\lambda$ by $\boldsymbol{t}_{n|\lambda}$ and $\boldsymbol{t}'_{n|\lambda}$ (resp., $\widetilde{\boldsymbol{t}}_{n|\lambda}$ and $\widetilde{\boldsymbol{t}}'_{n|\lambda}$). Then the following formulae are none other than valid versions of (3.7), (3.8), and (3.9) for $\boldsymbol{t}_{n|\lambda_n}$, $\boldsymbol{t}'_{n|\lambda_n}$, $\boldsymbol{T}_{n|\lambda_n}$, and their tilded counterparts, but all are developed making use of (2.36):

$$\boldsymbol{t}_{n|\lambda_n}(s)=\overline{\boldsymbol{T}_{n|\lambda_n}(s,\cdot)}=R^*_{\lambda_n}(T_n)\left(\boldsymbol{t}_n(s)\right)=\chi_n(s)R^*_{\lambda_n}(T_n)\left(\boldsymbol{t}(s)\right), \tag{3.11}$$

$$\boldsymbol{t}'_{n|\lambda_n}(t)=\boldsymbol{T}_{n|\lambda_n}(\cdot,t)=R_{\lambda_n}(T_n)\left(\boldsymbol{t}'_n(t)\right)=R_{\lambda_n}(T_n)P_n\left(\boldsymbol{t}'(t)\right), \tag{3.12}$$

$$\overline{\boldsymbol{T}_{n|\lambda_n}(s,t)}=\bar{\lambda}_n\left\langle\boldsymbol{t}_{n|\lambda_n}(s),P_n\left(\boldsymbol{t}'(t)\right)\right\rangle+\overline{\boldsymbol{T}_n(s,t)}, \tag{3.13}$$

$$\boldsymbol{T}_{n|\lambda_n}(s,t)=\lambda_n\chi_n(s)\left\langle\boldsymbol{t}'_{n|\lambda_n}(t),\boldsymbol{t}(s)\right\rangle+\boldsymbol{T}_n(s,t), \tag{3.14}$$

$$\begin{gathered}\widetilde{\boldsymbol{t}}_{n|\lambda_n}(s)=\overline{\widetilde{\boldsymbol{T}}_{n|\lambda_n}(s,\cdot)}=R^*_{\lambda_n}(\widetilde{T}_n)(\widetilde{\boldsymbol{t}}_n(s))=\chi_n(s)R^*_{\lambda_n}(\widetilde{T}_n)P_n\left(\boldsymbol{t}(s)\right),\\ \widetilde{\boldsymbol{t}}'_{n|\lambda_n}(t)=\widetilde{\boldsymbol{T}}_{n|\lambda_n}(\cdot,t)=R_{\lambda_n}(\widetilde{T}_n)(\widetilde{\boldsymbol{t}}'_n(t))=\chi_n(t)R_{\lambda_n}(\widetilde{T}_n)P_n\left(\boldsymbol{t}'(t)\right),\end{gathered} \tag{3.15}$$

$$\begin{gathered}\overline{\widetilde{\boldsymbol{T}}_{n|\lambda_n}(s,t)}=\bar{\lambda}_n\chi_n(t)\langle\widetilde{\boldsymbol{t}}_{n|\lambda_n}(s),P_n\left(\boldsymbol{t}'(t)\right)\rangle+\overline{\widetilde{\boldsymbol{T}}_n(s,t)},\\ \widetilde{\boldsymbol{T}}_{n|\lambda_n}(s,t)=\lambda_n\chi_n(s)\langle\widetilde{\boldsymbol{t}}'_{n|\lambda_n}(t),P_n\left(\boldsymbol{t}(s)\right)\rangle+\widetilde{\boldsymbol{T}}_n(s,t),\end{gathered} \tag{3.16}$$

where $s$, $t\in\mathbb{R}$ are arbitrary. It is readily seen from (2.27) and (3.14) that each $K^0$-kernel $\boldsymbol{T}_{n|\lambda}(s,t)$ has bounded $s$-support (namely, lying in $[-\tau_n,\tau_n]$), so the (in general, hardly verifiable) condition (3.3) of Definition 2 is automatically fulfilled with $\boldsymbol{T}_{n|\lambda}$ in the role of $\boldsymbol{T}_{|\lambda}$. Thus, in this role, $\boldsymbol{T}_{n|\lambda}$ is the only solution of the simultaneous integral equations (3.1) and (3.2) (with, of course, $\boldsymbol{T}$ replaced by $\boldsymbol{T}_n$) which is a $K^0$-kernel. The problem of explicitly finding that solution in terms of $\boldsymbol{T}_n$ is completely solved via the Fredholm-determinant method, as follows. For $\boldsymbol{T}_n$ a subkernel of $\boldsymbol{T}$, consider its Fredholm determinant $D_{\boldsymbol{T}_n}(\lambda)$ defined by the series

$$D_{\boldsymbol{T}_n}(\lambda):=1+\sum_{m=1}^{\infty}\frac{(-\lambda)^m}{m!}\int\cdots\int\boldsymbol{T}_n\begin{pmatrix}x_1&\dots&x_m\\x_1&\dots&x_m\end{pmatrix}dx_1\dots dx_m, \tag{3.17}$$

for every $\lambda\in\mathbb{C}$, and its first Fredholm minor $D_{\boldsymbol{T}_n}(s,t\mid\lambda)$ defined by the series

$$D_{\boldsymbol{T}_n}(s,t\mid\lambda)=\boldsymbol{T}_n(s,t)+\sum_{m=1}^{\infty}\frac{(-\lambda)^m}{m!}\int\cdots\int\boldsymbol{T}_n\begin{pmatrix}s&x_1&\dots&x_m\\t&x_1&\dots&x_m\end{pmatrix}dx_1\dots dx_m, \tag{3.18}$$

for all points $s$, $t\in\mathbb{R}$ and for every $\lambda\in\mathbb{C}$, where

$$\boldsymbol{T}_n\begin{pmatrix}x_1&\dots&x_\nu\\y_1&\dots&y_\nu\end{pmatrix}:=\det\begin{pmatrix}\boldsymbol{T}_n(x_1,y_1)&\dots&\boldsymbol{T}_n(x_1,y_\nu)\\\dots&\dots&\dots\\\boldsymbol{T}_n(x_\nu,y_1)&\dots&\boldsymbol{T}_n(x_\nu,y_\nu)\end{pmatrix}.$$

The next proposition can reliably be inferred from results of the Carleman-Mikhlin-Smithies theory for the Fredholm determinant and the first Fredholm minor for the Hilbert-Schmidt kernels with possibly unbounded support (see [3], [20], and [33]).

**Proposition 3.** *Let $\lambda \in \mathbb{C}$ be arbitrary but fixed. Then*

*1) the series of (3.17) is absolutely convergent in $\mathbb{C}$, and the series of (3.18) is absolutely convergent in $C(\mathbb{R}^2, \mathbb{C})$ and in $L^2(\mathbb{R}^2)$;*

*2) if $D_{\boldsymbol{T}_n}(\lambda) \neq 0$ then resolvent kernel for $\boldsymbol{T}_n$ at $\lambda$ exists and is the quotient of the first Fredholm minor and the Fredholm determinant:*

$$\boldsymbol{T}_{n|\lambda}(s,t) \equiv \frac{D_{\boldsymbol{T}_n}(s,t \mid \lambda)}{D_{\boldsymbol{T}_n}(\lambda)}; \tag{3.19}$$

*3) if $D_{\boldsymbol{T}_n}(\lambda) = 0$, then the resolvent kernel for $\boldsymbol{T}_n$ does not exist at $\lambda$.*

For each $n \in \mathbb{N}$, therefore, the characteristic set $\Lambda(T_n)(= \Lambda_p(T_n))$ is composed of all the zeros of the entire function $D_{\boldsymbol{T}_n}(\lambda)$, and is at most a denumerable set of complex numbers clustering at infinity. For $\widetilde{\boldsymbol{T}}_{n|\lambda}$, a formula like (3.19) is built up in the same way but replacing $\boldsymbol{T}_n$ by $\widetilde{\boldsymbol{T}}_n$. Since $\widetilde{T}_n^m = T_n^m P_n$ for $m \in \mathbb{N}$, the $m$-th iterants of $\widetilde{\boldsymbol{T}}_n$ and $\boldsymbol{T}_n$ (see (2.25)) stand therefore in a similar relation to each other, namely: $\widetilde{\boldsymbol{T}}_n^{[m]}(s,t) = \chi_n(t)\boldsymbol{T}_n^{[m]}(s,t)$ for all $s$, $t \in \mathbb{R}$. Then it follows from the rules for calculating the coefficients of powers of $\lambda$ in the Fredholm series (3.17) and (3.18) that $D_{\widetilde{\boldsymbol{T}}_n}(\lambda) \equiv D_{\boldsymbol{T}_n}(\lambda)$ and $D_{\widetilde{\boldsymbol{T}}_n}(s,t \mid \lambda) \equiv \chi_n(t) D_{\boldsymbol{T}_n}(s,t \mid \lambda)$. Hence, for each $n$,

$$\Lambda(T_n) = \Lambda(\widetilde{T}_n), \quad \Pi(T_n) = \Pi(\widetilde{T}_n), \tag{3.20}$$

$$\widetilde{\boldsymbol{T}}_{n|\lambda}(s,t) = \frac{D_{\widetilde{\boldsymbol{T}}_n}(s,t \mid \lambda)}{D_{\widetilde{\boldsymbol{T}}_n}(\lambda)} = \chi_n(t)\boldsymbol{T}_{n|\lambda}(s,t), \tag{3.21}$$

$$\widetilde{\boldsymbol{t}}_{n|\lambda}(s) = P_n\left(\boldsymbol{t}_{n|\lambda}(s)\right), \quad \widetilde{\boldsymbol{t}}'_{n|\lambda}(t) = \chi_n(t)\boldsymbol{t}'_{n|\lambda}(t).$$

Formal use of the triangle inequalities then gives the pointwise inequalities

$$\begin{aligned} |\widetilde{\boldsymbol{T}}_{n|\lambda}(s,t) - \boldsymbol{R}(s,t)| &\leqslant \chi_n(t)\left|\boldsymbol{T}_{n|\lambda}(s,t) - \boldsymbol{R}(s,t)\right| + \widehat{\chi}_n(t)|\boldsymbol{R}(s,t)|, \\ \|\widetilde{\boldsymbol{t}}_{n|\lambda}(s) - \boldsymbol{a}(s)\| &\leqslant \left\|P_n(\boldsymbol{t}_{n|\lambda}(s) - \boldsymbol{a}(s))\right\| + \|(I - P_n)(\boldsymbol{a}(s))\|, \\ \|\widetilde{\boldsymbol{t}}'_{n|\lambda}(t) - \boldsymbol{b}(t)\| &\leqslant \chi_n(t)\left\|\boldsymbol{t}'_{n|\lambda}(t) - \boldsymbol{b}(t)\right\| + \widehat{\chi}_n(t)\|\boldsymbol{b}(t)\|, \end{aligned} \tag{3.22}$$

where $\boldsymbol{R}\colon \mathbb{R}^2 \to \mathbb{C}$ and $\boldsymbol{a}, \boldsymbol{b}\colon \mathbb{R} \to L^2$ are arbitrary functions.

In summarizing the above listed material, we can say that, when being stated for $T = T_n$ (resp., $\widetilde{T}_n$), Korotkov's problem $(*)$ in Section 1 has as a fully satisfactory answer the Fredholm-type formula (3.19) (resp., (3.21)). Common to $\boldsymbol{T}_{n|\lambda}(s,t)$ and $\widetilde{\boldsymbol{T}}_{n|\lambda}(s,t)$ is an analytic feature of being written as a ratio of explicit entire functions in $\lambda$ for fixed values of $s$ and $t$, but in the present paper the feature will be of no essential use in proving the approximation results for the resolvent kernels $\boldsymbol{T}_{|\lambda}$ associated to the original kernel $\boldsymbol{T}$.

## 4. The main results

4.1. **On relative compactnesses.** Given an arbitrary bounded (possibly convergent) sequence $\{\lambda_n\}_{n=1}^{\infty}$ of complex numbers satisfying $\lambda_n \in \Pi(T_n)$ for each $n$, the $C(\mathbb{R}^2, \mathbb{C})$-valued sequence of the resolvent kernels

$$\{\boldsymbol{T}_{n|\lambda_n}\}_{n=1}^{\infty} \tag{4.1}$$

(all of whose terms, recall, are known analytically in terms of the original kernel $\boldsymbol{T}$ via the Fredholm formulae (3.17)-(3.19)) and the $C(\mathbb{R}, L^2)$-valued sequences,

$$\{\boldsymbol{t}_{n|\lambda_n}\}_{n=1}^{\infty} \quad \text{and} \quad \{\boldsymbol{t}'_{n|\lambda_n}\}_{n=1}^{\infty}, \tag{4.2}$$

of the respective Carleman resolvent functions are not known to be relatively compact in general. The following result gives a simple sufficient condition for the sequence (4.1) to be so.

**Theorem 3.** *Suppose that* $1 \in \nabla_{\mathfrak{b}}(\{\lambda_n T_n\})$*. Then the sequence* (4.1) *is relatively compact in* $C(\mathbb{R}^2, \mathbb{C})$*, and hence each sequence of positive integers contains a subsequence* $\{n_k\}_{k=1}^{\infty}$ *along which* $\boldsymbol{T}_{n_k|\lambda_{n_k}}$ *converges in the norm of* $C(\mathbb{R}^2, \mathbb{C})$*.*

*Proof.* Let $C > 0$ be that constant for which $|\lambda_n| \leqslant C$ for all $n$ in $\mathbb{N}$, and let (see (2.10) and Remark 1) $M = M(1) + 1$ be that constant for which $\|R_{\lambda_n}(T_n)\|(= \|R_1(\lambda_n T_n)\|) \leqslant M$ for all $n$ in $\mathbb{N}$ ($N(1) < 1$ for simplicity). From Lemma 2 applied with $R^*_{\lambda_n}(T_n)$ in place of $B_n$ and its proof it follows via (3.11) that the estimates

$$\|\boldsymbol{t}_{n|\lambda_n}(x)\| \leqslant M\|\boldsymbol{\tau}\|_{C(\mathbb{R},\mathbb{R})} \quad \text{and} \quad \|\boldsymbol{t}'_{n|\lambda_n}(x)\| \leqslant M\|\boldsymbol{\tau}'\|_{C(\mathbb{R},\mathbb{R})} \tag{4.3}$$

hold for all $x$ in $\mathbb{R}$ and all $n$ in $\mathbb{N}$ (see (2.41)), and that the subsets

$$\mathbf{t} = \bigcup_{n=1}^{\infty} \boldsymbol{t}_{n|\lambda_n} \subset C(\mathbb{R}, L^2) \quad \text{and} \quad \mathbf{t}' = \bigcup_{n=1}^{\infty} \boldsymbol{t}'_{n|\lambda_n} \subset C(\mathbb{R}, L^2) \tag{4.4}$$

are equicontinuous on $\mathbb{R} \cup \{\infty\}$, and, as seen from the proof of (2.42), their equicontinuity is *uniform* on $\mathbb{R}$, in the sense that for a fixed $\varepsilon$, the value of $\delta_x(\varepsilon)$ which is being used in (2.1) can be found to be the same for every $x$ in $\mathbb{R}$. Since in addition, $\boldsymbol{t}'$ and $\boldsymbol{T}$ are uniformly continuous vanishing at infinity, it follows that given any $\varepsilon > 0$, a $\delta > 0$ can be found such that for $x$, $y$, $u$ and $v$ in $\mathbb{R}$,

$$\max\left\{\sup_{n\in\mathbb{N}} \|\boldsymbol{t}_{n|\lambda_n}(x) - \boldsymbol{t}_{n|\lambda_n}(y)\|, \|\boldsymbol{t}'(u) - \boldsymbol{t}'(v)\|, |\boldsymbol{T}(x,u) - \boldsymbol{T}(y,v)|\right\} < \varepsilon$$

whenever $\max\{|x-y|, |u-v|\} < \delta$ or/and whenever both $\min\{|x|, |u|\} > \dfrac{1}{\delta}$ and $\boldsymbol{\tau}(y) = \boldsymbol{\tau}'(v) = |\boldsymbol{T}(y,v)| = 0$, whence, for all the $x$, $y$, $u$ and $v$ values so selected, the inequality

$$|\boldsymbol{T}_{n|\lambda_n}(x,u) - \boldsymbol{T}_{n|\lambda_n}(y,v)| \leqslant C\varepsilon\|\boldsymbol{\tau}'\|_{C(\mathbb{R},\mathbb{R})} + CM\|\boldsymbol{\tau}\|_{C(\mathbb{R},\mathbb{R})}\varepsilon + \varepsilon \tag{4.5}$$

follows, by means of the chain of inequalities that uses equation (3.13), takes into account property (2.29) and holds for all $x$, $y$, $u$, and $v$ in $\mathbb{R}$ and all $n$ in $\mathbb{N}$:

$$\begin{aligned}
&|\boldsymbol{T}_{n|\lambda_n}(x,u) - \boldsymbol{T}_{n|\lambda_n}(y,v)| \\
&\leqslant |\lambda_n| \left|\left\langle \boldsymbol{t}_{n|\lambda_n}(x), P_n\left(\boldsymbol{t}'(u)\right)\right\rangle - \left\langle \boldsymbol{t}_{n|\lambda_n}(y), P_n\left(\boldsymbol{t}'(v)\right)\right\rangle\right| + |\boldsymbol{T}_n(x,u) - \boldsymbol{T}_n(y,v)| \\
&\leqslant |\lambda_n| \left|\left\langle \boldsymbol{t}_{n|\lambda_n}(x) - \boldsymbol{t}_{n|\lambda_n}(y), P_n\left(\boldsymbol{t}'(u)\right)\right\rangle\right| + |\lambda_n| \left|\left\langle \boldsymbol{t}_{n|\lambda_n}(y), P_n\left(\boldsymbol{t}'(u) - \boldsymbol{t}'(v)\right)\right\rangle\right| \\
&\quad + |\boldsymbol{T}_n(x,u) - \boldsymbol{T}_n(y,v)| \\
&\leqslant |\lambda_n| \|\boldsymbol{t}_{n|\lambda_n}(x) - \boldsymbol{t}_{n|\lambda_n}(y)\| \boldsymbol{\tau}'(u) + |\lambda_n| \|\boldsymbol{t}_{n|\lambda_n}(y)\| \|\boldsymbol{t}'(u) - \boldsymbol{t}'(v)\| \\
&\quad + \begin{cases} \max\{|\boldsymbol{T}(x,u)|, |\boldsymbol{T}(y,v)|\} & \text{if } |x-y| = |\tau_n - x| + |\tau_n - y|, \\ \max\{|\boldsymbol{T}(x,u)|, |\boldsymbol{T}(y,v)|\} & \text{if } |x-y| = |\tau_n + x| + |\tau_n + y|, \\ |\boldsymbol{T}(x,u) - \boldsymbol{T}(y,v)| & \text{otherwise.} \end{cases}
\end{aligned}$$

Observe that this chain, when written for fixed $y$ and $v$ satisfying $\boldsymbol{\tau}(y) = \boldsymbol{\tau}'(v) = |\boldsymbol{T}(y,v)| = 0$, also leads (via (4.3)) to the following estimate, valid for all $x$ and $u$ in $\mathbb{R}$ and all $n$ in $\mathbb{N}$:

$$|\boldsymbol{T}_{n|\lambda_n}(x,u)| \leqslant 3CM\|\boldsymbol{\tau}\|_{C(\mathbb{R},\mathbb{R})}\|\boldsymbol{\tau}'\|_{C(\mathbb{R},\mathbb{R})} + \|\boldsymbol{T}\|_{C(\mathbb{R}^2,\mathbb{C})}. \tag{4.6}$$

Now, by (4.5), the subset

$$\mathfrak{T} = \bigcup_{n=1}^{\infty} \boldsymbol{T}_{n|\lambda_n} \subset C(\mathbb{R}^2, \mathbb{C}) \tag{4.7}$$

is equicontinuous on $\mathbb{R}^2 \cup \{\infty\}$, i.e., it does have property (A2) of Theorem 1, for $k = 2$ and $B = \mathbb{C}$. By (4.6), the subset $\mathfrak{T}$ in (4.7) does have property (A1) of that theorem, also with $\mathbb{R}^2$ and $\mathbb{C}$ in place of $\mathbb{R}^k$ and $B$, respectively. Linking properties (A1) and (A2) yields, via the same Theorem 1, the relative compactness of $\mathfrak{T}$ in $C(\mathbb{R}^2, \mathbb{C})$. The theorem is proved. □

The theorem just proved gives some support to the following two questions: if the sequence (4.1) converges in $C(\mathbb{R}^2, \mathbb{C})$, possibly along some subsequence, say $\{\boldsymbol{T}_{n_k|\lambda_{n_k}}\}$, toward a function $\boldsymbol{R}$ say, (1) whether $\lambda = \lim_{k\to\infty} \lambda_{n_k}$ is necessarily a regular value for $T$, and (2) if $\lambda$ does turn out to belong to $\Pi(T)$, whether $\boldsymbol{R} = \boldsymbol{T}_{|\lambda}$. With the advantage of being more closely related to the $L^2$-topology of the underlying space, similar questions can and should be raised for the $C(\mathbb{R}, L^2)$-valued sequences of Carleman functions in (4.2). But, in order for these sequences to be relatively compact in $C(\mathbb{R}, L^2)$ (even if, as in the previous proof, the subsets $\mathbf{t}$, $\mathbf{t'} \subset C(\mathbb{R}, L^2)$ are equicontinuous on $\mathbb{R} \cup \{\infty\}$ and bounded) the following two conditions of pointwise type have to be satisfied, respectively:

(A) the $L^2$-valued sequence

$$\{\boldsymbol{t}_{n|\lambda_n}(s)\}_{n=1}^{\infty} \tag{4.8}$$

is relatively compact in $L^2$ for each fixed $s \in \mathbb{R}$;

(B) the $L^2$-valued sequence

$$\left\{\boldsymbol{t}'_{n|\lambda_n}(t)\right\}_{n=1}^{\infty} \tag{4.9}$$

is relatively compact in $L^2$ for each fixed $t \in \mathbb{R}$.

An investigation into the consequences of conditions (A) and (B) is the chief purpose of the remainder of the present subsection. Here is the main result, containing also some possible answers to the above questions, but under some weakened conditions on the convergence of the sequences under consideration.

**Theorem 4.** *Suppose that the regular values $\lambda_n \in \Pi(T_n)$, involved in (4.1), (4.8), and (4.9), are chosen so that*

$$\lim_{n\to\infty} \lambda_n = \lambda \tag{4.10}$$

*for some non-zero $\lambda \in \mathbb{C}$. Then the following assertions are true:*

(1) *If* (A) *holds, then $\lambda \notin \Lambda_p(T)$.*

(2) *If* (B) *holds, then $\lambda \notin \Lambda_r(T)$.*

(3) *If both* (A) *and* (B) *hold, then*

(a) *the sequence* (4.1) *has a pointwise limit, $\boldsymbol{R}$:*

$$\boldsymbol{R}(s,t) = \lim_{n\to\infty} \boldsymbol{T}_{n|\lambda_n}(s,t) \quad \textit{for every } (s,t) \in \mathbb{R}^2; \tag{4.11}$$

(b) *the integral equation of the second kind,*

$$f(s) - \lambda \int \boldsymbol{T}(s,t) f(t)\, dt = g(s) \quad \textit{for almost every } s \in \mathbb{R}, \tag{4.12}$$

*is solvable in $f \in L^2$ if and only if the function $g$ of $L^2$ in the right side of the equation and the pointwise limit function $\boldsymbol{R}$ in* (4.11) *stand in the following two relations to each other:*

$$\int \left| \int \boldsymbol{R}(s,t) g(t)\, dt \right|^2 ds < \infty, \tag{4.13}$$

$$\int \boldsymbol{T}(s,x)\int \boldsymbol{R}(x,t)g(t)\,dt\,dx = \int g(t)\int \boldsymbol{T}(s,x)\boldsymbol{R}(x,t)\,dx\,dt \quad \textit{for every } s\in\mathbb{R}, \tag{4.14}$$

*and in that case the equation* (4.12) *has a unique solution* $f \in L^2$ *given by the equation*

$$f(s) = g(s) + \lambda \int \boldsymbol{R}(s,t)g(t)\,dt \quad \textit{for almost every } s\in\mathbb{R}. \tag{4.15}$$

*Remark* 7. The solvability conditions (4.13)-(4.14) strongly resemble those which are labeled (b)-(c) in Theorem of [26], although the latter are expressed by means of generalized Fredholm minors.

*Proof.* (1) Under assumption (A), to each temporarily fixed $s \in \mathbb{R}$ may be assigned at least one function $\boldsymbol{a}(s) \in L^2$ obtainable as the $L^2$ limit of the sequence (4.8) where $n$ is restricted to an appropriate sequence $\{p_n\}_{n=1}^{\infty}$ of positive integers increasing to $\infty$. (Of course, it may happen that different choices of the sequence $\{p_n\}$ generate different functions $\boldsymbol{a}(s)$.) With $s$ kept fixed and with $n$ replaced by $p_n$, let equation (3.13) be rewritten in terms of Carleman functions and operators as the following equation in $L^2$:

$$\left(I - \bar{\lambda}_{p_n} T_{p_n}^{*}\right)\left(\boldsymbol{t}_{p_n|\lambda_{p_n}}(s)\right) = \boldsymbol{t}_{p_n}(s)$$

(compare also (3.11)). By letting $n$ tend to $\infty$, this equation becomes

$$\left(I - \bar{\lambda} T^{*}\right)\left(\boldsymbol{a}(s)\right) = \boldsymbol{t}(s), \tag{4.16}$$

due to (4.10), (2.31), and (2.35). Arbitrariness of $s\in\mathbb{R}$ then implies that $\boldsymbol{t}(\mathbb{R}) \subset \operatorname{Ran}\left(I - \bar{\lambda}T^{*}\right)$, whence it follows via (2.24) that $\overline{\operatorname{Ran}(T^{*})} \subset \overline{\operatorname{Ran}\left(I - \bar{\lambda}T^{*}\right)}$. In its turn, the latter inclusion implies that $\operatorname{Ker}(I - \lambda T) \subset \operatorname{Ker}(T)$, which is incompatible with $\lambda \in \Lambda_p(T)$, i.e., assertion (1) in the theorem is true.

(2) Supposing (B), to prove assertion (2) in the theorem, start by rewriting (in the same manner as above) equation (3.14) as

$$\left(I - \lambda_{r_n} T_{r_n}\right)\left(\boldsymbol{t}'_{r_n|\lambda_{r_n}}(t)\right) = \boldsymbol{t}'_{r_n}(t) \tag{4.17}$$

(compare also (3.12)), where $\{r_n\}_{n=1}^{\infty}$ is deliberately chosen to be a strictly increasing subsequence of positive integers, generally depending on $t$, such that the $L^2$ limit $\boldsymbol{b}(t) := \lim_{n\to\infty} \boldsymbol{t}'_{r_n|\lambda_{r_n}}(t)$ exists. Then, on account of (4.10), (2.31) and (2.35), letting $n$ in (4.17) go to infinity yields the equation

$$(I - \lambda T)\left(\boldsymbol{b}(t)\right) = \boldsymbol{t}'(t) \tag{4.18}$$

valid for all $t$ in $\mathbb{R}$. Together with (2.24), this equation implies that $\overline{\operatorname{Ran}(T)} \subset \overline{\operatorname{Ran}(I - \lambda T)}$, and therefore that $\operatorname{Ker}(I - \bar{\lambda}T^{*}) \subset \operatorname{Ker}(T^{*})$. The latter in turn implies that $\bar{\lambda} \notin \Lambda_p(T^{*})$. It then follows from (2.14) that $\lambda \notin \Lambda_r(T) \subset (\Lambda_p(T^{*}))^{\#}$, and assertion (2) is thus also true.

(3) If conditions (A) and (B) hold simultaneously, then assertions (1) and (2) together imply via (2.12) and (2.11) that either $\lambda \in \Lambda_c(T)$ or $\lambda \in \Pi(T)$. Hence, in any case, both the operators $I - \lambda T$ and $I - \bar{\lambda}T^{*}$ are one-to-one (see (2.14), (2.13), and (2.2)). Whenever $s \in \mathbb{R}$ (resp., $t \in \mathbb{R}$) is fixed, from (4.16) (resp., (4.18)) it then follows that every $L^2$-convergent subsequence of the sequence (4.8) (resp., (4.9)) has the same limit. Therefore, for each fixed $s$ and $t$, the relatively compact sequences (4.8) and (4.9) are themselves convergent, implying that there are $L^2$-valued functions $\boldsymbol{a}$ and $\boldsymbol{b}$, defined on all of $\mathbb{R}$, such that

$$\lim_{n\to\infty}\left\|\boldsymbol{a}(s) - \boldsymbol{t}_{n|\lambda_n}(s)\right\| = 0 \quad \text{and} \quad \lim_{n\to\infty}\left\|\boldsymbol{b}(t) - \boldsymbol{t}'_{n|\lambda_n}(t)\right\| = 0 \tag{4.19}$$

for every $s$, $t \in \mathbb{R}$. Passing to the limit as $n \to \infty$ is then possible in the right (and, as a consequence, in the left) sides of the equations (see (3.13) and (3.14))

$$\boldsymbol{T}_{n|\lambda_n}(s,t) = \lambda_n \langle P_n\left(\boldsymbol{t'}(t)\right), \boldsymbol{t}_{n|\lambda_n}(s)\rangle + \boldsymbol{T}_n(s,t),$$
$$\boldsymbol{T}_{n|\lambda_n}(s,t) = \lambda_n \chi_n(s)\overline{\langle \boldsymbol{t}(s), \boldsymbol{t'}_{n|\lambda_n}(t)\rangle} + \boldsymbol{T}_n(s,t)$$

thanks to (2.33), (2.30), and (4.10); it respectively results in the equations

$$\boldsymbol{R}(s,t) = \lambda \left\langle \boldsymbol{t'}(t), \boldsymbol{a}(s)\right\rangle + \boldsymbol{T}(s,t), \tag{4.20}$$
$$\boldsymbol{R}(s,t) = \lambda\overline{\langle \boldsymbol{t}(s), \boldsymbol{b}(t)\rangle} + \boldsymbol{T}(s,t) \tag{4.21}$$

valid for all $s$, $t \in \mathbb{R}$, where the function $\boldsymbol{R}\colon \mathbb{R}^2 \to \mathbb{C}$ can safely be regarded as being defined pointwise by

$$\boldsymbol{R}(s,t) = \lim_{n\to\infty} \boldsymbol{T}_{n|\lambda_n}(s,t) \tag{4.22}$$

at each $(s,t) \in \mathbb{R}^2$ (thereby proving (4.11)).

In the case when $\lambda$ is in $\Pi(T)$, one may simply compare formulae (3.7) with (4.16) and (4.18), and infer that $\boldsymbol{a} = \boldsymbol{t}_{|\lambda}$ and $\boldsymbol{b} = \boldsymbol{t'}_{|\lambda}$ in the $C(\mathbb{R}, L^2)$ sense, whence the equality $\boldsymbol{R} = \boldsymbol{T}_{|\lambda}$ in $C(\mathbb{R}^2, \mathbb{C})$ emerges by comparison of (3.8) and (3.9) with (4.20) and (4.21). But in the general case when $\lambda \in \Lambda_c(T) \cup \Pi(T)$, the proof of assertion (3)(b) requires somewhat deeper arguments, and is as follows. Comparing (4.19) with (4.22) and taking note of (3.11) and (3.12) yields that $\overline{\boldsymbol{R}(s,\cdot)} \in \boldsymbol{a}(s)$ and $\boldsymbol{R}(\cdot,t) \in \boldsymbol{b}(t)$ for every $s$, $t \in \mathbb{R}$ (inasmuch as pointwise and $L^2$ limits agree almost everywhere). Therefore, on writing the $L^2$ inner product as an integral, equations (4.20) and (4.21) go over into

$$\boldsymbol{R}(s,t) = \lambda \int \boldsymbol{R}(s,x)\boldsymbol{T}(x,t)\,dx + \boldsymbol{T}(s,t), \tag{4.23}$$
$$\boldsymbol{R}(s,t) = \lambda \int \boldsymbol{T}(s,x)\boldsymbol{R}(x,t)\,dx + \boldsymbol{T}(s,t), \tag{4.24}$$

respectively (compare with (3.2) and (3.1)). From the properties of $\boldsymbol{R}$ and $\boldsymbol{T}$ it follows that, after being multiplied by any $L^2$-function $h$, the equations just displayed are integrable with respect to $t$, so as to yield

$$\begin{aligned}\int \boldsymbol{R}(s,t)h(t)\,dt &= \lambda \int \left(\int \boldsymbol{R}(s,x)\boldsymbol{T}(x,t)\,dx\right) h(t)\,dt + \int \boldsymbol{T}(s,t)h(t)\,dt,\\ \int \boldsymbol{R}(s,t)h(t)\,dt &= \lambda \int \left(\int \boldsymbol{T}(s,x)\boldsymbol{R}(x,t)\,dx\right) h(t)\,dt + \int \boldsymbol{T}(s,t)h(t)\,dt,\end{aligned} \tag{4.25}$$

whence

$$\int \left(\int \boldsymbol{R}(s,x)\boldsymbol{T}(x,t)\,dx\right) h(t)\,dt = \int \left(\int \boldsymbol{T}(s,x)\boldsymbol{R}(x,t)\,dx\right) h(t)\,dt. \tag{4.26}$$

To prove the "only if" part of assertion (3)(b), let $f$ and $g$ be any two $L^2$-functions obeying (4.12) and suppose that the function $\boldsymbol{R}$ is as in (4.22) above and therefore satisfies both (4.23) and (4.24). Then it holds that

$$\begin{aligned}\int \boldsymbol{R}(s,t)g(t)\,dt &= \int \boldsymbol{R}(s,t)\left(f(t) - \lambda \int \boldsymbol{T}(t,x)f(x)\,dx\right) dt\\ &= \int \left(\boldsymbol{R}(s,t) - \lambda \int \boldsymbol{R}(s,x)\boldsymbol{T}(x,t)\,dx\right) f(t)\,dt\\ &= \int \boldsymbol{T}(s,t)f(t)\,dt = \frac{f(s)-g(s)}{\lambda},\end{aligned} \tag{4.27}$$

where the very last equality is in the almost everywhere sense, and where the reversal of the order of integration in deducing the second line from the first one may be proved in a way similar to (2.26), as follows:

$$\begin{aligned}\int \boldsymbol{R}(s,t)\int \boldsymbol{T}(t,x)f(x)\,dx\,dt &= \langle Tf,\boldsymbol{a}(s)\rangle = \langle f,T^{*}(\boldsymbol{a}(s))\rangle \\ &= \int f(t)\overline{\int \boldsymbol{T'}(t,x)\overline{\boldsymbol{R}(s,x)}\,dx}\,dt = \int\left(\int \boldsymbol{R}(s,x)\boldsymbol{T}(x,t)\,dx\right)f(t)\,dt.\end{aligned} \tag{4.28}$$

From (4.27) the asserted property (4.13) follows. Furthermore, from (4.27) it follows that

$$\begin{aligned}\int \boldsymbol{T}(s,x)\left(\int \boldsymbol{R}(x,t)g(t)\,dt\right)\,dx &= \int \boldsymbol{T}(s,x)\left(\int \boldsymbol{T}(x,t)f(t)\,dt\right)\,dx \\ &= \int \boldsymbol{T}^{[2]}(s,t)f(t)\,dt.\end{aligned} \tag{4.29}$$

On the other hand, with the aid of the subsequent use of (4.26), (4.12), the first equation in (4.25) with $h(t)=\int \boldsymbol{T}(t,x)f(x)\,dx$, and (4.28), one gets

$$\begin{aligned}\int\left(\int \boldsymbol{T}(s,x)\boldsymbol{R}(x,t)\,dx\right)g(t)\,dt &= \int\left(\int \boldsymbol{R}(s,x)\boldsymbol{T}(x,t)\,dx\right)g(t)\,dt \\ &= \int\left(\int \boldsymbol{R}(s,x)\boldsymbol{T}(x,t)\,dx\right)\left(f(t)-\lambda\int \boldsymbol{T}(t,x)f(x)\,dx\right)\,dt \\ &= \int\left(\int \boldsymbol{R}(s,x)\boldsymbol{T}(x,t)\,dx\right)f(t)\,dt \\ &\quad -\lambda\int\left(\int \boldsymbol{R}(s,x)\boldsymbol{T}(x,t)\,dx\right)\left(\int \boldsymbol{T}(t,x)f(x)\,dx\right)\,dt \\ &= \int\left(\int \boldsymbol{R}(s,x)\boldsymbol{T}(x,t)\,dx\right)f(t)\,dt \\ &\quad -\int \boldsymbol{R}(s,t)\left(\int \boldsymbol{T}(t,x)f(x)\,dx\right)\,dt + \int \boldsymbol{T}^{[2]}(s,t)f(t)\,dt \\ &= \int \boldsymbol{T}^{[2]}(s,t)f(t)\,dt.\end{aligned}$$

Whence property (4.14) follows via comparison with (4.29).

To prove the "if" part of assertion (3)(b), let both conditions (4.13) and (4.14) (with $\boldsymbol{R}$ given by (4.22)) be satisfied by a certain fixed function $g$ of $L^2$, and let the function $f$ be of the form (4.15). It is to be proved that $f$ (uniquely) solves (4.12). Substitute, therefore, the expression for $f$ into (4.12), to get, applying (4.14) and (4.24), that

$$\begin{aligned}f(s)-\lambda\int \boldsymbol{T}(s,t)f(t)\,dt &= g(s)-\lambda\int \boldsymbol{T}(s,t)g(t)\,dt \\ &\quad +\lambda\int \boldsymbol{R}(s,t)g(t)\,dt-\lambda^{2}\int \boldsymbol{T}(s,x)\left(\int \boldsymbol{R}(x,t)g(t)\,dt\right)\,dx \\ &= g(s)-\lambda\int \boldsymbol{T}(s,t)g(t)\,dt+\lambda\int\left(\boldsymbol{R}(s,t)-\lambda\int \boldsymbol{T}(s,x)\boldsymbol{R}(x,t)\,dx\right)g(t)\,dt \\ &= g(s)-\lambda\int \boldsymbol{T}(s,t)g(t)\,dt+\lambda\int \boldsymbol{T}(s,t)g(t)\,dt = g(s)\end{aligned}$$

almost everywhere in $\mathbb{R}$.

To prove the uniqueness, assume that $v \in L^2$ is also a solution of equation (4.12), i.e., assume that

$$v(s) - \lambda \int \boldsymbol{T}(s,t)v(t)dt = g(s) \quad \text{for almost every } s \text{ in } \mathbb{R}.$$

Then, using (4.15), (4.28), and (4.23), one obtains

$$\begin{aligned} f(s) &= g(s) + \lambda \int \boldsymbol{R}(s,t)g(t)\,dt \\ &= v(s) - \lambda \int \boldsymbol{T}(s,t)v(t)dt + \lambda \int \boldsymbol{R}(s,t)\left(v(t) - \lambda \int \boldsymbol{T}(t,x)v(x)\,dx\right)\,dt \\ &= v(s) - \lambda \int \boldsymbol{T}(s,t)v(t)dt + \lambda \int \boldsymbol{R}(s,t)v(t)\,dt - \lambda^2 \int \boldsymbol{R}(s,t) \int \boldsymbol{T}(t,x)v(x)dx\,dt \\ &= v(s) - \lambda \int \left(\boldsymbol{T}(s,t) - \boldsymbol{R}(s,t)\right)v(t)\,dt - \lambda^2 \int \left(\int \boldsymbol{R}(s,x)\boldsymbol{T}(x,t)\,dx\right)v(t)\,dt \\ &= v(s) - \lambda \int \left(\boldsymbol{T}(s,t) - \boldsymbol{R}(s,t) + \lambda \int \boldsymbol{R}(s,x)\boldsymbol{T}(x,t)\,dx\right)v(t)\,dt = v(s) \end{aligned}$$

almost everywhere in $\mathbb{R}$. The proof of the theorem is complete. □

**Corollary 1.** *In the conditions of Theorem 4, if the original kernel $\boldsymbol{T}$ induces a quasi-compact operator $T \in \mathfrak{R}(L^2)$, then condition* (A) *or* (B) *holds if and only if $\lambda \in \Pi(T)$, and in that case the sequences of* (4.2) *converge in $C(\mathbb{R}, L^2)$ to the resolvent Carleman functions $\boldsymbol{t}_{|\lambda}$ and $\boldsymbol{t}'_{|\lambda}$, respectively, and the sequence* (4.1) *converges in $C(\mathbb{R}^2, \mathbb{C})$ to the resolvent kernel $\boldsymbol{T}_{|\lambda}$.*

*Proof.* If the resolvent Carleman function sequence (4.8) (resp., (4.9)) is relatively compact in $L^2$ for each $s \in \mathbb{R}$ (resp., $t \in \mathbb{R}$), then, by Theorem 4, $\lambda \notin \Lambda_p(T)$ (resp., $\bar{\lambda} \notin \Lambda_p(T^*)$). This implies, on account of the quasi-compactness of $T$ and (2.15), that in either of these cases, $\lambda \in \Pi(T)$. Since, by the compactness of $T^m$ for some positive integer $m$,

$$\|(T - T_n)T_n^m\| \leqslant \|T(T_n^m - T^m)\| + \|T^{m+1} - T_n^{m+1}\| \to 0 \quad \text{as } n \to \infty, \tag{4.30}$$

the "if" part of the corollary follows from Theorems 5 and 6 of Subsection 4.2, putting there $\beta_n = \frac{\lambda_n - \lambda}{\lambda \lambda_n}$. The corollary is proved. □

4.2. **Parameter range conditions for uniform approximations.** The question we deal with in this subsection is in a sense converse to those we have touched upon in the previous one, and is as follows: given that $\lambda$ is as defined in Theorem 4 and is a (non-zero) regular value for $T$, what further connections between $\{\lambda_n\}$ and $\lambda$ guarantee the existence, in suitable uniform senses, of the following limit relations involving the whole sequences:

$$\boldsymbol{t}_{|\lambda} = \lim_{n\to\infty} \boldsymbol{t}_{n|\lambda_n}, \quad \boldsymbol{t}'_{|\lambda} = \lim_{n\to\infty} \boldsymbol{t}'_{n|\lambda_n}, \quad \boldsymbol{T}_{|\lambda} = \lim_{n\to\infty} \boldsymbol{T}_{n|\lambda_n}?$$

In the theorem which follows, we characterize one such connection by means of Kato's regions like $\nabla_{\mathfrak{s}}(\cdot)$, defined at the end of Subsection 2.3.

**Theorem 5.** *Let $\{\beta_n\}_{n=1}^{\infty}$ be an arbitrary sequence of complex numbers satisfying*

$$\lim_{n\to\infty} \beta_n = 0, \tag{4.31}$$

*and define $\lambda_n(\lambda) := \lambda(1 - \beta_n \lambda)^{-1}$, so that one can consider that $\lambda_n(\lambda) \to \lambda$ when $n \to \infty$ for each fixed $\lambda \in \mathbb{C}$. Then $\varnothing \neq \nabla_{\mathfrak{s}}(\{\beta_n I + T_n\}) \subseteq \nabla_{\mathfrak{s}}(\{\beta_n I + \widetilde{T}_n\}) \subseteq \Pi(T)$, and the following limits hold true whenever $\mathfrak{K}$ (resp., $\widetilde{\mathfrak{K}}$) is a compact subset of $\nabla_{\mathfrak{s}}(\{\beta_n I + T_n\}$ (resp., $\nabla_{\mathfrak{s}}(\{\beta_n I + \widetilde{T}_n\})$):*

$$\lim_{n\to\infty} \sup_{\lambda \in \widetilde{\mathfrak{K}}} \left\| \boldsymbol{t}'_{n|\lambda_n(\lambda)} - \boldsymbol{t}'_{|\lambda} \right\|_{C(\mathbb{R}, L^2)} = 0, \tag{4.32}$$

$$\lim_{n\to\infty}\sup_{\lambda\in\mathfrak{K}}\left\|\boldsymbol{t}_{n|\lambda_n(\lambda)}-\boldsymbol{t}_{|\lambda}\right\|_{C(\mathbb{R},L^2)}=0, \tag{4.33}$$

$$\lim_{n\to\infty}\sup_{\lambda\in\widetilde{\mathfrak{K}}}\left\|\boldsymbol{T}_{n|\lambda_n(\lambda)}-\boldsymbol{T}_{|\lambda}\right\|_{C(\mathbb{R}^2,\mathbb{C})}=0. \tag{4.34}$$

*Remark* 8. The formulae (4.32)-(4.34) strengthen the convergence properties for the sequences (4.1) and (4.2) that are presented in [28] by six analogous limit relations, numbered (63)-(68), each holding uniformly with respect to only one of the variables for all fixed values of the other. The proof, however, is almost the same, being also based upon Lemma 1.

*Proof.* Let us begin by collecting (mainly from [14]) some preparatory results, to be numbered below from (4.36) to (4.46). To shorten notation, write $A_n := \beta_n I + T_n$, $\widetilde{A}_n := \beta_n I + \widetilde{T}_n$. Choose a (non-zero) regular value $\zeta \in \Pi(T)$ satisfying $|\zeta|\,\|T\| < 1$, and hence satisfying for some $N(\zeta) > 0$ the inequality

$$|\zeta|\,\|A_n\| \leqslant |\zeta|\left(\max_{n>N(\zeta)}|\beta_n| + \|T\|\right) < 1 \quad \text{for all } n > N(\zeta). \tag{4.35}$$

Then $\zeta$ does belong to $\nabla_{\mathfrak{b}}(\{A_n\})$, because

$$\left\|A_{n|\zeta}\right\| \leqslant \frac{\|A_n\|}{1-|\zeta|\,\|A_n\|} \leqslant M(\zeta) = \frac{\max\limits_{n>N(\zeta)}\|A_n\|}{1-|\zeta|\left(\max\limits_{n>N(\zeta)}|\beta_n|+\|T\|\right)} \quad \text{for all } n > N(\zeta)$$

(compare (2.10)). The result is that the intersection of $\nabla_{\mathfrak{b}}(\{A_n\})$ and $\Pi(T)$ is non-void. Similarly it can be shown that $\nabla_{\mathfrak{b}}(\{\widetilde{A}_n\}) \cap \Pi(T) \neq \varnothing$. Therefore, since, because of (4.31) and (2.31), the sequences $\{A_n\}$ and $\{\widetilde{A}_n\}$ both converge to $T$ in the strong operator topology, it follows by the criterion for generalized strong convergence (see [14, Theorem VIII-1.5]) that

$$\nabla_{\mathfrak{s}}(\{A_n\}) = \nabla_{\mathfrak{b}}(\{A_n\}) \cap \Pi(T), \quad \nabla_{\mathfrak{s}}(\{\widetilde{A}_n\}) = \nabla_{\mathfrak{b}}(\{\widetilde{A}_n\}) \cap \Pi(T), \tag{4.36}$$

$$\lim_{n\to\infty}\left\|\left(A_{n|\lambda}-T_{|\lambda}\right)f\right\| = 0 \quad \text{for all } \lambda \in \nabla_{\mathfrak{s}}(\{A_n\}) \text{ and } f \in L^2,$$

$$\lim_{n\to\infty}\|(\widetilde{A}_{n|\lambda}-T_{|\lambda})f\| = 0 \quad \text{for all } \lambda \in \nabla_{\mathfrak{s}}(\{\widetilde{A}_n\}) \text{ and } f \in L^2. \tag{4.37}$$

Further, given a $\lambda \in \nabla_{\mathfrak{b}}(\{A_n\}) \cup \nabla_{\mathfrak{b}}(\{\widetilde{A}_n\})$, the following formulae hold for sufficiently large $n$:

$$\begin{aligned}
R_\lambda(A_n) &= \frac{1}{1-\beta_n\lambda}R_{\lambda_n(\lambda)}(T_n), \quad R_\lambda(\widetilde{A}_n) = \frac{1}{1-\beta_n\lambda}R_{\lambda_n(\lambda)}(\widetilde{T}_n),\\
R_{\lambda_n(\lambda)}(T_n) &= (1-\beta_n\lambda)\left(I+\lambda A_{n|\lambda}\right) = I + \lambda\beta_n I + \lambda A_{n|\lambda} + \lambda^2\beta_n A_{n|\lambda},\\
A_{n|\lambda} &= \left(\frac{1}{1-\beta_n\lambda}\right)^2 T_{n|\lambda_n(\lambda)} + \frac{\beta_n}{1-\beta_n\lambda}I,\\
\widetilde{A}_{n|\lambda} &= \left(\frac{1}{1-\beta_n\lambda}\right)^2 \widetilde{T}_{n|\lambda_n(\lambda)} + \frac{\beta_n}{1-\beta_n\lambda}I.
\end{aligned} \tag{4.38}$$

These are obtained by a purely formal calculation, and use that fact that $\Pi(\widetilde{T}_n) = \Pi(T_n)$ for each fixed $n \in \mathbb{N}$ (see (3.20)). The equations in the last two lines combine to give, using (3.21),

$$A_{n|\lambda}P_n = \widetilde{A}_{n|\lambda} + \frac{\beta_n}{1-\beta_n\lambda}(I-P_n). \tag{4.39}$$

This implies in particular that $\|\widetilde{A}_{n|\lambda}\| \leqslant \left\|A_{n|\lambda}\right\| + |\beta_n||1-\beta_n\lambda|^{-1}$, whence (4.31) leads to the inclusion relation $\nabla_{\mathfrak{b}}(\{A_n\}) \subseteq \nabla_{\mathfrak{b}}(\{\widetilde{A}_n\})$, from which it follows via (4.36) that $\varnothing \neq \nabla_{\mathfrak{s}}(\{A_n\}) \subseteq \nabla_{\mathfrak{s}}(\{\widetilde{A}_n\}) \subseteq \Pi(T)$, as asserted.

In what follows, let $L$ denote a relatively compact set in $L^2$, and let $\widetilde{\mathfrak{K}}$ denote a compact (closed and bounded) subset of $\nabla_{\mathfrak{s}}(\{\widetilde{A}_n\})$. Then, according to Theorem VIII-1.1 in [14] there exists a positive constant $M(\widetilde{\mathfrak{K}})$ such that

$$\sup_{\lambda\in\widetilde{\mathfrak{K}}}\|\widetilde{A}_{n|\lambda}\| \leqslant M(\widetilde{\mathfrak{K}}) \quad \text{for all sufficiently large } n, \tag{4.40}$$

so the convergence in (4.37) is uniform over both $\lambda\in\widetilde{\mathfrak{K}}$ and $f\in L$:

$$\lim_{n\to\infty}\sup_{\substack{\lambda\in\widetilde{\mathfrak{K}}\\ f\in L}}\|(\widetilde{A}_{n|\lambda}-T_{|\lambda})f\| = 0. \tag{4.41}$$

Indeed, the following relations hold:

$$\begin{aligned}
&\lim_{n\to\infty}\sup_{\substack{\lambda\in\widetilde{\mathfrak{K}}\\ f\in L}}\|(\widetilde{A}_{n|\lambda}-T_{|\lambda})f\| \\
&= \lim_{n\to\infty}\sup_{\substack{\lambda\in\widetilde{\mathfrak{K}}\\ f\in L}}\|(I+\lambda\widetilde{A}_{n|\lambda})(T-\widetilde{A}_n)R_\lambda(T)f\| && \text{by (2.6)}\\
&\leqslant \sup_{\lambda\in\widetilde{\mathfrak{K}}}(1+|\lambda|M(\widetilde{\mathfrak{K}}))\lim_{n\to\infty}\sup_{\substack{\lambda\in\widetilde{\mathfrak{K}}\\ f\in L}}\left\|(T-\widetilde{A}_n)R_\lambda(T)f\right\| && \text{by (4.40)}\\
&= 0 && \text{by Lemma 1,}
\end{aligned}$$

inasmuch as the set

$$\bigcup_{\lambda\in\widetilde{\mathfrak{K}};\, f\in L} R_\lambda(T)f$$

is relatively compact in $L^2$; this is due to (2.7), the compactness of $\mathfrak{K}$ and the relative compactness of $L$.

Now use (4.40), (4.41), and the observation from (4.31) that

$$\sup_{\lambda\in\widetilde{\mathfrak{K}}}\left|\frac{\beta_n}{1-\beta_n\lambda}\right| \leqslant \frac{|\beta_n|}{1-|\beta_n|\sup\limits_{\lambda\in\widetilde{\mathfrak{K}}}|\lambda|} \to 0 \quad \text{as } n\to\infty, \tag{4.42}$$

to infer, via the connecting formula (4.39), that

$$\sup_{\lambda\in\widetilde{\mathfrak{K}}}\left\|A_{n|\lambda}P_n\right\| < M(\widetilde{\mathfrak{K}})+1 \quad \text{for all sufficiently large } n, \tag{4.43}$$

$$\lim_{n\to\infty}\sup_{\substack{\lambda\in\widetilde{\mathfrak{K}}\\ f\in L}}\left\|\left(A_{n|\lambda}P_n-T_{|\lambda}\right)f\right\| = 0. \tag{4.44}$$

Throughout what follows let $\mathfrak{K}$ denote a compact subset of $\nabla_{\mathfrak{s}}(\{A_n\})$. Then Theorem VIII-1.1 in [14], this time applied to the operator sequence $\{A_n\}$, yields the conclusion that there exists a positive constant $M(\mathfrak{K})$ such that

$$\sup_{\lambda\in\mathfrak{K}}\left\|A_{n|\lambda}\right\| \leqslant M(\mathfrak{K}) \quad \text{for all sufficiently large } n. \tag{4.45}$$

Hence there holds

$$\lim_{n\to\infty}\sup_{\substack{\lambda\in\mathfrak{K}\\ f\in L}}\left\|\left(A_{n|\lambda}-T_{|\lambda}\right)^* f\right\| = 0. \tag{4.46}$$

Indeed, the following relations hold:

$$\begin{aligned}
&\lim_{n\to\infty}\sup_{\substack{\lambda\in\mathfrak{K}\\ f\in L}}\left\|\left(A_{n|\lambda}-T_{|\lambda}\right)^{*}f\right\| \\
&=\lim_{n\to\infty}\sup_{\substack{\lambda\in\mathfrak{K}\\ f\in L}}\left\|\left(I+\bar{\lambda}\left(A_{n|\lambda}\right)^{*}\right)\left(T-A_n\right)^{*}R^{*}_{\lambda}(T)f\right\| && \text{by (2.6)}\\
&\leqslant \sup_{\lambda\in\mathfrak{K}}(1+|\lambda|M(\mathfrak{K}))\lim_{n\to\infty}\sup_{\substack{\lambda\in\mathfrak{K}\\ f\in L}}\left\|\left(T-A_n\right)^{*}R^{*}_{\lambda}(T)f\right\| && \text{by (4.45)}\\
&=0 && \text{by Lemma 1,}
\end{aligned}$$

inasmuch as $(A_n)^* \to T^*$ strongly as $n \to \infty$ (by (2.31), (4.31)) and the set

$$\bigcup_{\lambda\in\mathfrak{K};\, f\in L} R^{*}_{\lambda}(T)f$$

is relatively compact in $L^2$ (by the compactness of $\mathfrak{K}$, (2.7) and the relative compactness of $L$). Similarly, it can be proved that

$$\lim_{n\to\infty}\sup_{\substack{\lambda\in\widetilde{\mathfrak{K}}\\ f\in L}}\left\|P_n\left(A_{n|\lambda}P_n-T_{|\lambda}\right)^{*}f\right\|=0. \tag{4.47}$$

With these preparations, we are ready to establish that the limit formulae (4.32)-(4.34) all hold. To this effect, use formulae (3.12), (3.11), (3.7), (4.38), and then the triangle inequality to write

$$\begin{aligned}
&\sup_{\lambda\in\widetilde{\mathfrak{K}}}\left\|\boldsymbol{t}'_{n|\lambda_n(\lambda)}-\boldsymbol{t}'_{|\lambda}\right\|_{C(\mathbb{R},L^2)}\\
&\quad=\sup_{\lambda\in\widetilde{\mathfrak{K}}}\sup_{t\in\mathbb{R}}\left\|\left(P_n-I+\lambda\beta_nP_n+\lambda(A_{n|\lambda}P_n-T_{|\lambda})+\lambda^2\beta_nA_{n|\lambda}P_n\right)\left(\boldsymbol{t}'(t)\right)\right\|\\
&\qquad\leqslant\sup_{t\in\mathbb{R}}\left\|(P_n-I)\left(\boldsymbol{t}'(t)\right)\right\|+|\beta_n|\sup_{\lambda\in\widetilde{\mathfrak{K}}}|\lambda|\cdot\sup_{t\in\mathbb{R}}\left\|P_n\left(\boldsymbol{t}'(t)\right)\right\|\\
&\qquad+\sup_{\substack{\lambda\in\widetilde{\mathfrak{K}}\\ t\in\mathbb{R}}}\left(|\lambda|\left\|\left(A_{n|\lambda}P_n-T_{|\lambda}\right)\left(\boldsymbol{t}'(t)\right)\right\|\right)+|\beta_n|\|\boldsymbol{\tau}'\|_{C(\mathbb{R},\mathbb{R})}\sup_{\lambda\in\widetilde{\mathfrak{K}}}\left(|\lambda|^2\|A_{n|\lambda}P_n\|\right),
\end{aligned}\tag{4.48}$$

$$\begin{aligned}
&\sup_{\lambda\in\mathfrak{K}}\left\|\boldsymbol{t}_{n|\lambda_n(\lambda)}-\boldsymbol{t}_{|\lambda}\right\|_{C(\mathbb{R},L^2)}\\
&\quad=\sup_{\lambda\in\mathfrak{K}}\sup_{s\in\mathbb{R}}\left\|\left(\chi_n(s)\left(I+\lambda\beta_nI+\lambda A_{n|\lambda}+\lambda^2\beta_nA_{n|\lambda}\right)-I-\lambda T_{|\lambda}\right)^{*}\left(\boldsymbol{t}(s)\right)\right\|\\
&\qquad\leqslant|\beta_n|\|\chi_n\boldsymbol{\tau}\|_{C(\mathbb{R},\mathbb{R})}\sup_{\lambda\in\mathfrak{K}}|\lambda|+\|\widehat{\chi}_n\boldsymbol{\tau}\|_{C(\mathbb{R},\mathbb{R})}\\
&\qquad+\sup_{\substack{\lambda\in\mathfrak{K}\\ s\in\mathbb{R}}}\left(\chi_n(s)\,|\lambda|\left\|\left(A_{n|\lambda}-T_{|\lambda}\right)^{*}\left(\boldsymbol{t}(s)\right)\right\|\right)+\|\widehat{\chi}_n\boldsymbol{\tau}\|_{C(\mathbb{R},\mathbb{R})}\sup_{\lambda\in\mathfrak{K}}\left(|\lambda|\left\|T_{|\lambda}\right\|\right)\\
&\qquad+|\beta_n|\|\chi_n\boldsymbol{\tau}\|_{C(\mathbb{R},\mathbb{R})}\sup_{\lambda\in\mathfrak{K}}(|\lambda|^2\left\|A_{n|\lambda}\right\|).
\end{aligned}\tag{4.49}$$

Because of (2.30), (4.44) with $L=\boldsymbol{t}'(\mathbb{R})$, (4.43), (4.31), (2.38), and of the boundedness of the set $\widetilde{\mathfrak{K}}$, each summand on the right-hand side of (4.48) tends as $n\to\infty$ to zero. This proves (4.32). To see that formula (4.33) also holds, apply (4.46) with $L=\boldsymbol{t}(\mathbb{R})$, (4.45), (2.22), (4.31), and (2.38) to the right-hand-side terms in (4.49), hereby proving their convergence to zero.

Now use equations (3.9), (3.14), and the triangle and the Cauchy-Schwarz inequality to also write

$$\begin{aligned}
&\sup_{\lambda\in\widetilde{\mathfrak{K}}}\left\|\boldsymbol{T}_{n|\lambda_n(\lambda)}-\boldsymbol{T}_{|\lambda}\right\|_{C(\mathbb{R}^2,\mathbb{C})}= \\
&\quad=\sup_{\lambda\in\widetilde{\mathfrak{K}}}\sup_{(s,t)\in\mathbb{R}^2}\left|\lambda_n(\lambda)\chi_n(s)\left\langle \boldsymbol{t}'_{n|\lambda_n(\lambda)}(t),\boldsymbol{t}(s)\right\rangle-\lambda\left\langle \boldsymbol{t}'_{|\lambda}(t),\boldsymbol{t}(s)\right\rangle+\boldsymbol{T}_n(s,t)-\boldsymbol{T}(s,t)\right| \\
&\qquad\leqslant \sup_{\lambda\in\widetilde{\mathfrak{K}}}\sup_{(s,t)\in\mathbb{R}^2}\left(\chi_n(s)\,|\lambda|\left|\left\langle \boldsymbol{t}'_{n|\lambda_n(\lambda)}(t)-\boldsymbol{t}'_{|\lambda}(t),\boldsymbol{t}(s)\right\rangle\right|\right) \\
&\qquad+\sup_{\lambda\in\widetilde{\mathfrak{K}}}\sup_{(s,t)\in\mathbb{R}^2}\left(\chi_n(s)\,|\lambda_n(\lambda)-\lambda|\left|\left\langle \boldsymbol{t}'_{n|\lambda_n(\lambda)}(t),\boldsymbol{t}(s)\right\rangle\right|\right) \\
&\qquad+\sup_{\lambda\in\widetilde{\mathfrak{K}}}\sup_{(s,t)\in\mathbb{R}^2}\left(\widehat{\chi}_n(s)\,|\lambda|\left|\left\langle \boldsymbol{t}'_{|\lambda}(t),\boldsymbol{t}(s)\right\rangle\right|\right)+\sup_{(s,t)\in\mathbb{R}^2}\left|\boldsymbol{T}_n(s,t)-\boldsymbol{T}(s,t)\right| \\
&\qquad\leqslant \|\chi_n\boldsymbol{\tau}\|_{C(\mathbb{R},\mathbb{R})}\sup_{\lambda\in\widetilde{\mathfrak{K}}}|\lambda|\cdot\sup_{\lambda\in\widetilde{\mathfrak{K}}}\left\|\boldsymbol{t}'_{n|\lambda_n(\lambda)}-\boldsymbol{t}'_{|\lambda}\right\|_{C(\mathbb{R},L^2)} \\
&\qquad+\|\chi_n\boldsymbol{\tau}\|_{C(\mathbb{R},\mathbb{R})}\sup_{\lambda\in\widetilde{\mathfrak{K}}}\left(\left|\frac{\lambda^2\beta_n}{1-\beta_n\lambda}\right|\left\|\boldsymbol{t}'_{n|\lambda_n(\lambda)}\right\|_{C(\mathbb{R},L^2)}\right) \\
&\qquad+\|\widehat{\chi}_n\boldsymbol{\tau}\|_{C(\mathbb{R},\mathbb{R})}\|\boldsymbol{\tau}'\|_{C(\mathbb{R},\mathbb{R})}\sup_{\lambda\in\widetilde{\mathfrak{K}}}\left(|\lambda|\,\|R_\lambda(T)\|\right)+\|\boldsymbol{T}_n-\boldsymbol{T}\|_{C(\mathbb{R}^2,\mathbb{C})}.
\end{aligned}$$

The above-established relations (2.33), (4.32), (4.42), and (2.22) together imply that the four terms on the extreme right side here all converge to zero as $n\to\infty$, which proves (4.34). The theorem is proved. □

*Remark* 9. Applying the respective results of Theorem 5 in conjunction with the inequalities (compare (3.22))

$$\begin{aligned}
&\sup_{\lambda\in\widetilde{\mathfrak{K}}}\left\|\widetilde{\boldsymbol{t}}'_{n|\lambda_n(\lambda)}-\boldsymbol{t}'_{|\lambda}\right\|_{C(\mathbb{R},L^2)}\leqslant\sup_{\lambda\in\widetilde{\mathfrak{K}}}\left\|\boldsymbol{t}'_{n|\lambda_n(\lambda)}-\boldsymbol{t}'_{|\lambda}\right\|_{C(\mathbb{R},L^2)}+\|\widehat{\chi}_n\boldsymbol{\tau}'\|_{C(\mathbb{R},\mathbb{R})}\sup_{\lambda\in\widetilde{\mathfrak{K}}}\|R_\lambda(T)\|, \\
&\sup_{\lambda\in\widetilde{\mathfrak{K}}}\left\|\widetilde{\boldsymbol{T}}_{n|\lambda_n(\lambda)}-\boldsymbol{T}_{|\lambda}\right\|_{C(\mathbb{R}^2,\mathbb{C})}\leqslant\sup_{\lambda\in\widetilde{\mathfrak{K}}}\left\|\boldsymbol{T}_{n|\lambda_n(\lambda)}-\boldsymbol{T}_{|\lambda}\right\|_{C(\mathbb{R}^2,\mathbb{C})}+\sup_{(s,t)\in\mathbb{R}^2}\left|\widehat{\chi}_n(t)\boldsymbol{T}(s,t)\right| \\
&\qquad+\|\widehat{\chi}_n\boldsymbol{\tau}'\|_{C(\mathbb{R},\mathbb{R})}\|\boldsymbol{\tau}\|_{C(\mathbb{R},\mathbb{R})}\cdot\sup_{\lambda\in\widetilde{\mathfrak{K}}}|\lambda|
\end{aligned}$$

yields that the limits (4.32) and (4.34) both remain valid upon replacing $\boldsymbol{t}'_{n|\lambda_n(\lambda)}$ and $\boldsymbol{T}_{n|\lambda_n(\lambda)}$ by $\widetilde{\boldsymbol{t}}'_{n|\lambda_n(\lambda)}$ and $\widetilde{\boldsymbol{T}}_{n|\lambda_n(\lambda)}$, respectively. In turn, the limit (4.33) continues to hold with $\boldsymbol{t}_{n|\lambda_n(\lambda)}$ and $\mathfrak{K}$ $(\subset\nabla_{\mathfrak{s}}(\{\beta_nI+T_n\}))$ replaced respectively by $\widetilde{\boldsymbol{t}}_{n|\lambda_n(\lambda)}$ and $\widetilde{\mathfrak{K}}$ $(\subset\nabla_{\mathfrak{s}}(\{\beta_nI+\widetilde{T}_n\}))$, and to prove this use can be made of the inequality (compare (4.49))

$$\begin{aligned}
&\sup_{\lambda\in\widetilde{\mathfrak{K}}}\left\|\widetilde{\boldsymbol{t}}_{n|\lambda_n(\lambda)}-\boldsymbol{t}_{|\lambda}\right\|_{C(\mathbb{R},L^2)}\leqslant|\beta_n|\|\chi_n\boldsymbol{\tau}\|_{C(\mathbb{R},\mathbb{R})}\sup_{\lambda\in\widetilde{\mathfrak{K}}}|\lambda|+\|\widehat{\chi}_n\boldsymbol{\tau}\|_{C(\mathbb{R},\mathbb{R})} \\
&\qquad+\sup_{\substack{\lambda\in\widetilde{\mathfrak{K}}\\ s\in\mathbb{R}}}\left(\chi_n(s)\,|\lambda|\left\|P_n\left(A_{n|\lambda}-T_{|\lambda}\right)^*(\boldsymbol{t}(s))\right\|\right)+\|\widehat{\chi}_n\boldsymbol{\tau}\|_{C(\mathbb{R},\mathbb{R})}\sup_{\lambda\in\widetilde{\mathfrak{K}}}\left(|\lambda|\left\|T_{|\lambda}\right\|\right) \\
&\qquad+|\beta_n|\|\chi_n\boldsymbol{\tau}\|_{C(\mathbb{R},\mathbb{R})}\sup_{\lambda\in\widetilde{\mathfrak{K}}}(|\lambda|^2\left\|A_{n|\lambda}P_n\right\|)+\sup_{\substack{\lambda\in\widetilde{\mathfrak{K}}\\ s\in\mathbb{R}}}\left\|(I-P_n)\left(\boldsymbol{t}_{|\lambda}(s)\right)\right\|,
\end{aligned}$$

of (4.47) applied with $L = \boldsymbol{t}(\mathbb{R})$, of (4.43), and of the relative compactness in $L^2$ of the set

$$\bigcup_{\lambda\in\widetilde{\mathfrak{K}};\, s\in\mathbb{R}} \boldsymbol{t}_{|\lambda}(s).$$

4.3. **Remarks on approximations at the whole set of regular values.** In connection with Theorem 5 the following natural question can be asked: in what cases are the sets $\overset{\circ}{\Pi}(T) := \Pi(T)\setminus\{0\}$ and $\nabla_{\mathfrak{s}}(\{\frac{\lambda_n-\lambda}{\lambda\lambda_n}I+T_n\})$ coincident? One answer to this question is given in the following theorem.

**Theorem 6.** *Let $m$ be a fixed positive integer, and for each $n\in\mathbb{N}$, let $J_n$ be that nuclear operator on $L^2$ which is defined by the $K^0$-kernel*

$$\boldsymbol{J}_n(s,t) = \widehat{\chi}_n(s)\int_{\mathbb{I}_n} \boldsymbol{T}(s,x)\boldsymbol{T}_n^{[m]}(x,t)\,dx \quad (s,\, t\in\mathbb{R}). \tag{4.50}$$

*Then, if*

$$\|J_n\| \to 0 \quad \text{as } n\to\infty, \tag{4.51}$$

*for any complex sequence $\{\beta_n\}$ converging to zero there holds*

$$\nabla_{\mathfrak{s}}(\{\beta_n I + T_n\}) = \overset{\circ}{\Pi}(T) \subset \nabla_{\mathfrak{b}}(\{\beta_n I + T_n\}). \tag{4.52}$$

*Proof.* Rewrite condition (4.51) by (4.50) and the multipliability property of $K^0$-kernels (see (2.26)) as

$$\|(T-T_n)T_n^m\| \to 0 \quad \text{as } n\to\infty. \tag{4.53}$$

Continue to denote $A_n := \beta_n I + T_n$ as in the previous proof. Let $\lambda$ be a fixed non-zero regular value for $T$. A straightforward calculation yields the equation

$$\begin{aligned}\left((I-\lambda T)\sum_{k=0}^{m-1}\lambda^k A_n^k + \lambda^m A_n^m\right)(I-\lambda A_n)\\ = (I-\lambda T)\left(I+\lambda^{m+1}R_\lambda(T)(T-A_n)A_n^m\right).\end{aligned} \tag{4.54}$$

Expanding binomially $(\beta_n I + T_n)^m$ and utilizing conditions (4.31) and (4.53) gives

$$\begin{aligned}\|(T-A_n)A_n^m\| &= \|(T-\beta_n I - T_n)\,(\beta_n I+T_n)^m\|\\ &\leqslant \|(T-T_n)T_n^m\| + \|\beta_n T_n^m\|\\ &+ \|(T-\beta_n I-T_n)\|\sum_{k=1}^{m}\binom{m}{k}\left|\beta_n^k\right|\left\|T_n^{m-k}\right\| \to 0 \quad \text{as } n\to\infty,\end{aligned}$$

so $|\lambda|^{m+1}\|R_\lambda(T)(T-A_n)A_n^m\| < \frac{1}{2}$ for all $n$ sufficiently large. Note that, for such $n$, the right-hand side of equation (4.54) does represent an invertible operator on $L^2$. This makes the last factor

$$I-\lambda A_n = (1-\beta_n\lambda)\left(I-\tfrac{\lambda}{1-\beta_n\lambda}T_n\right) \tag{4.55}$$

on the left-hand side one-to-one and so invertible, as $T_n$ is compact. Hence, for such $n$, $\frac{\lambda}{1-\beta_n\lambda}\in\Pi(T_n)$, $\lambda\in\Pi(A_n)$, and

$$\left\|A_{n|\lambda}\right\| = \frac{1}{|\lambda|}\|R_\lambda(A_n)-I\|$$

$$= \frac{1}{|\lambda|}\left\|\left(I+\lambda^{m+1}R_\lambda(T)(T-A_n)A_n^m\right)^{-1}R_\lambda(T)\left((1-\lambda T)\sum_{k=0}^{m-1}(\lambda A_n)^k + (\lambda A_n)^m\right) - I\right\|$$

$$\leqslant \frac{1}{|\lambda|} \frac{\|R_\lambda(T)\|\,(1+|\lambda|\,\|T\|)\sum\limits_{k=0}^{m}|\lambda|^k\,\|A_n\|^k}{1-|\lambda|^{m+1}\,\|R_\lambda(T)(T-A_n)A_n^m\|}+\frac{1}{|\lambda|}$$

$$\leqslant M(\lambda) := \frac{2}{|\lambda|}\,\|R_\lambda(T)\|\,(1+|\lambda|\,\|T\|)\sum_{k=0}^{m}|\lambda|^k\,(\max_{n\in\mathbb{N}}|\beta_n|+\|T\|)^k+\frac{1}{|\lambda|},$$

where in the second equality use has been made of equation (4.54). Thus (see (2.10)), $\lambda \in \nabla_{\mathfrak{b}}(\{A_n\})$, and (4.52) now follows by (4.36). The theorem is proved. □

*Remark* 10. Condition like (4.53) also occurs in the theory of spectral approximation of linear operators. As proved, e.g., in [1, Theorem 2.2], an operator sequence $\{T_n\}$ converging strongly to an operator $T$ and satisfying (4.53) is necessarily a *strongly stable approximation* (following the terminology of F. Chatelin [5]) of the operator $T$ at an isolated spectral point in the boundary of the spectrum of $T$ (see also [21]).

*Remark* 11. Observe by (2.28) that (4.53) implies that

$$\|(T-\widetilde{T}_n)\widetilde{T}_n^m\| = \|(T-T_n)T_n^m P_n\| \to 0 \quad \text{as } n\to\infty, \tag{4.56}$$

which means that the $K^0$-kernels $\widetilde{\boldsymbol{J}}_n$,

$$\widetilde{\boldsymbol{J}}_n(s,t) = \widehat{\chi}_n(s)\int_{\mathbb{I}_n}\boldsymbol{T}(s,x)\widetilde{\boldsymbol{T}}_n^{[m]}(x,t)\,dx,$$

do induce on $L^2$ the (nuclear) operators $\widetilde{J}_n$, $\widetilde{J}_n=(T-\widetilde{T}_n)\widetilde{T}_n^m$, that also have their operator norm going to zero as $n$ goes to infinity. The same conclusion as (4.52) can be made from (4.56) for the sequence $\{\widetilde{T}_n\}$, it reads as follows:

$$\nabla_{\mathfrak{s}}(\{\beta_n I+\widetilde{T}_n\}) = \mathring{\Pi}(T) \subset \nabla_{\mathfrak{b}}(\{\beta_n I+\widetilde{T}_n\})$$

for every complex null sequence $\{\beta_n\}$. Consequently, under condition (4.53),

$$\nabla_{\mathfrak{s}}(\{\beta_n I+T_n\}) = \nabla_{\mathfrak{s}}(\{\beta_n I+\widetilde{T}_n\}) = \mathring{\Pi}(T) \tag{4.57}$$

for every complex null sequence $\{\beta_n\}$. In particular, if the nuclear operators $(I-P_n)TP_n$, which have as kernels the $K^0$-kernels $\widehat{\chi}_n(s)\boldsymbol{T}(s,t)\chi_n(t)$, converge to zero operator in the operator norm as $n\to\infty$, then both conditions (4.53) and (4.56) automatically hold for every fixed $m$ in $\mathbb{N}$. This holds, e.g., if $\boldsymbol{T}$ has the same support as that of $\boldsymbol{K}$ in Example 1, and this can happen even if $T$ is not a quasi-compact operator.

Now we are in a position to get the following reinforcement of the "if" part of Corollary 1.

**Theorem 7.** *If the original kernel $\boldsymbol{T}$ or some of its iterant $\boldsymbol{T}^{[m]}$ ($m \geqslant 2$) defines a compact operator on $L^2$, then the following limits hold true whenever $\mathfrak{K}$ is a compact subset of $\mathring{\Pi}(T)$:*

$$\begin{aligned}
&\lim_{n\to\infty}\sup_{\lambda\in\mathfrak{K}}\left\|\boldsymbol{t}'_{n|\lambda}-\boldsymbol{t}'_{|\lambda}\right\|_{C(\mathbb{R},L^2)}=0, && \lim_{n\to\infty}\sup_{\lambda\in\mathfrak{K}}\left\|\widetilde{\boldsymbol{t}}'_{n|\lambda}-\boldsymbol{t}'_{|\lambda}\right\|_{C(\mathbb{R},L^2)}=0,\\
&\lim_{n\to\infty}\sup_{\lambda\in\mathfrak{K}}\|\boldsymbol{t}_{n|\lambda}-\boldsymbol{t}_{|\lambda}\|_{C(\mathbb{R},L^2)}=0, && \lim_{n\to\infty}\sup_{\lambda\in\mathfrak{K}}\|\widetilde{\boldsymbol{t}}_{n|\lambda}-\boldsymbol{t}_{|\lambda}\|_{C(\mathbb{R},L^2)}=0,\\
&\lim_{n\to\infty}\sup_{\lambda\in\mathfrak{K}}\left\|\boldsymbol{T}_{n|\lambda}-\boldsymbol{T}_{|\lambda}\right\|_{C(\mathbb{R}^2,\mathbb{C})}=0, && \lim_{n\to\infty}\sup_{\lambda\in\mathfrak{K}}\|\widetilde{\boldsymbol{T}}_{n|\lambda}-\boldsymbol{T}_{|\lambda}\|_{C(\mathbb{R}^2,\mathbb{C})}=0.
\end{aligned} \tag{4.58}$$

*Proof.* Observe from (4.30) that under the stated hypotheses on the kernel $\boldsymbol{T}$, its operator $T$ is subject to the condition (4.51). With $\beta_n$ all taken equal to zero, apply Theorems 6, 5 and Remarks 9, 11 to conclude that the limits in (4.58) all hold. The theorem is proved. □

*Remark* 12. It is interesting to consider whether Theorem 7 could be strengthened to the statement that if the original kernel $\boldsymbol{T}$ defines an asymptotically quasi-compact compact operator $T$ on $L^2$, then the limits in (4.58) all hold.

4.4. **On finding characteristic values and spectral functions.** Throughout this subsection we shall assume that the original $K^0$-kernel $\boldsymbol{T}$ is Hermitian (i.e., $\boldsymbol{T}(s,t) = \overline{\boldsymbol{T}(t,s)}$ for all $s$, $t \in \mathbb{R}$), so that the integral operator $T \in \mathfrak{R}(L^2)$ it induces is self-adjoint ($T = T^*$) and hence has a strongly countably additive resolution of the identity $E_T \colon \Omega \to \mathfrak{P}(L^2)$, defined on the $\sigma$-algebra $\Omega$ of all Borel subsets of $\mathbb{R}$, such that $T = \int_{\mathbb{R}} z\, E_T(dz)$. Let $\Omega_0$ be the collection of all Borel sets $\omega \in \Omega$ whose closures $\overline{\omega}$ in $\mathbb{R}$ do not contain the point $\lambda = 0$. If $\chi_\omega$ is the characteristic function of the set $\omega \in \Omega_0$, then $\chi_\omega(\lambda) = \lambda v_\omega(\lambda)$, where $v_\omega(\lambda) := \chi_\omega(\lambda)/\lambda$. By the multiplicative property of the spectral measure $E_T$, it follows that $E_T(\omega) = \chi_\omega(T) = Tv_\omega(T) = v_\omega(T)T$, where $v_\omega(T) = \int_\omega \frac{1}{\lambda} E_T(d\lambda) \in \mathfrak{R}(L^2)$. Then Proposition 1 implies that the orthogonal projection $E_T(\omega)$ is a bi-Carleman operator such that the $C(\mathbb{R}, L^2)$-functions

$$\boldsymbol{e}_{|\omega}(\cdot) = v_\omega(T)(\boldsymbol{t}(\cdot)) \quad \text{and} \quad \boldsymbol{e}'_{|\omega}(\cdot) = v_\omega(T)(\boldsymbol{t}'(\cdot)) = \boldsymbol{e}_{|\omega}(\cdot) \tag{4.59}$$

are just its associated Carleman functions. Since $E_T(\omega) = (E_T(\omega))^2$, the bi-Carleman kernel $\boldsymbol{E}_{|\omega}$ associated to $E_T(\omega)$ can be computed pointwise as follows by means of (2.25):

$$\boldsymbol{E}_{|\omega}(s,t) = \boldsymbol{E}^{[2]}_{|\omega}(s,t) = \langle \boldsymbol{e}_{|\omega}(t), \boldsymbol{e}_{|\omega}(s)\rangle; \tag{4.60}$$

hence it is in $C(\mathbb{R}^2, \mathbb{C})$. We are therefore led to frame the following definition.

**Definition 3.** A set function $\boldsymbol{E}_{|\,\cdot}\colon \Omega_0 \to C(\mathbb{R}^2, \mathbb{C})$ whose value on each $\omega \in \Omega_0$ is the inducing $K^0$-kernel $\boldsymbol{E}_{|\omega}$ of the orthogonal projection $E_T(\omega)$ is called a *spectral function* for $\boldsymbol{T}$.

When $\omega$ is a one-point subset of $\mathbb{R}$, the following theorem gives a technique of using the double sequences $\{\widetilde{\boldsymbol{T}}_{m|\lambda_n}\}_{m,n=1}^\infty$ and $\{\widetilde{\boldsymbol{t}}_{m|\lambda_n}\}_{m,n=1}^\infty$ to get uniform approximations for $\boldsymbol{E}_{|\omega}$ and $\boldsymbol{e}_{|\omega}$ and to recognize whether the point contained in $\omega$ is the reciprocal of a characteristic value for $T$.

**Theorem 8.** *Let $\lambda$ be a non-zero real number, and let $\{\mu_n\}_{n=1}^\infty$ be a decreasing sequence of positive real numbers such that*

$$\lim_{n\to\infty} \mu_n = 0. \tag{4.61}$$

*Then a strictly increasing sequence $\{m_n\}_{n=1}^\infty$ of positive integers can be found in terms of the original $K^0$-kernel $\boldsymbol{T}$, along which, as $n \to \infty$,*

$$-\imath\mu_n \widetilde{\boldsymbol{t}}_{m_n|\lambda+\imath\mu_n} \to \boldsymbol{e}_{|\{\frac{1}{\lambda}\}} \quad \text{in } C(\mathbb{R}, L^2), \tag{4.62}$$

$$\imath\mu_n \widetilde{\boldsymbol{T}}_{m_n|\lambda+\imath\mu_n} \to \boldsymbol{E}_{|\{\frac{1}{\lambda}\}} \quad \text{in } C(\mathbb{R}^2, \mathbb{C}), \tag{4.63}$$

*and along which, therefore, the limit $\lim\limits_{n\to\infty} \mu_n \|\widetilde{\boldsymbol{T}}_{m_n|\lambda+\imath\mu_n}\|_{C(\mathbb{R}^2,\mathbb{C})}$ exists, which, according to (2.16), is equal to (resp., greater than) zero if and only if $\lambda \in \Lambda_c(T) \cup \Pi(T)$ (resp., $\Lambda_p(T)$). E.g., as such a sequence of indices one can always take a strictly increasing sequence $\{m_n\}_{n=1}^\infty$ of positive integers satisfying*

$$\lim_{n\to\infty} \sup_{s\in\mathbb{I}_{m_n}} \int_{\mathbb{I}_{m_n}} \left| \int_{\mathbb{I}_{m_n}} \boldsymbol{T}(t,\xi)\overline{\widetilde{\boldsymbol{T}}_{m_n|\lambda+\imath\mu_n}(s,\xi)}\, d\xi \right|^2 dt = 0. \tag{4.64}$$

*Proof.* For each $n \in \mathbb{N}$, let $\lambda_n = \lambda + \imath\mu_n$, where as before $\imath$ denotes the imaginary unit. First observe that, with $n$ fixed, $(T - \widetilde{T}_m)R^*_{\lambda_n}(\widetilde{T}_m)P_m \to 0$ strongly as $m \to \infty$; the convergence can be proved as follows, using (2.5), (2.9), (2.31), and (2.30): for $f \in L^2$,

$$\begin{aligned}\|(T-\widetilde{T}_m)R^*_{\lambda_n}(\widetilde{T}_m)P_m f\| &= \|\left(I-\bar{\lambda}_n T\right)R^*_{\lambda_n}(T)(T-\widetilde{T}_m)R^*_{\lambda_n}(\widetilde{T}_m)P_m f\| \\ &= \|\left(I-\bar{\lambda}_n T\right)R^*_{\lambda_n}(\widetilde{T}_m)(T-\widetilde{T}_m)R^*_{\lambda_n}(T)P_m f\| \\ &\leqslant \|I-\bar{\lambda}_n T\|\|R^*_{\lambda_n}(\widetilde{T}_m)\|\|(T-\widetilde{T}_m)R^*_{\lambda_n}(T)P_m f\| \\ \leqslant \frac{|\lambda_n|}{\mu_n}&(1+|\lambda_n|\,\|T\|)\|(T-\widetilde{T}_m)R^*_{\lambda_n}(T)P_m f\| \to 0 \quad \text{as } m\to\infty.\end{aligned}$$

Now since, by (3.15) and (2.38),

$$\begin{aligned}\sup_{s\in\mathbb{I}_m}\|(T-\widetilde{T}_m)(\widetilde{\boldsymbol{t}}_{m|\lambda_n}(s))\| &= \sup_{s\in\mathbb{I}_m}\left(\chi_m(s)\|(T-\widetilde{T}_m)R^*_{\lambda_n}(\widetilde{T}_m)P_m(\boldsymbol{t}(s))\|\right) \\ &\leqslant \sup_{s\in\mathbb{R}}\|(T-\widetilde{T}_m)R^*_{\lambda_n}(\widetilde{T}_m)P_m(\boldsymbol{t}(s))\| \to 0 \quad \text{as } m\to\infty,\end{aligned}$$

it follows that given an index $n\in\mathbb{N}$, an index $m_n\in\mathbb{N}$ can be found such that, purely in terms of the known $K^0$-kernels (see (3.21)),

$$\begin{aligned}\frac{1}{n} \geqslant \varepsilon_n := \sup_{s\in\mathbb{I}_{m_n}}&\sqrt{\int_{\mathbb{I}_{m_n}}\left|\int_{\mathbb{I}_{m_n}}\boldsymbol{T}(t,x)\overline{\widetilde{\boldsymbol{T}}_{m_n|\lambda_n}(s,x)}\,dx\right|^2 dt} \\ \left(= \sup_{s\in\mathbb{I}_{m_n}}\right.&\left.\|(T-\widetilde{T}_{m_n})(\widetilde{\boldsymbol{t}}_{m_n|\lambda_n}(s))\|\right)\end{aligned} \tag{4.65}$$

(compare (4.64)). Clearly the sequence $\{m_n\}_{n=1}^{\infty}$ may be chosen to be strictly increasing.

Recall (from, e.g., p. 1366 of [7]) that the operator sequence $\{-\imath\frac{\mu_n}{\lambda_n}R^*_{\lambda_n}(T)\}_{n=1}^{\infty}$ converges to

$$E_T\left(\{\tfrac{1}{\lambda}\}\right)\ \left(=\lambda T E_T\left(\{\tfrac{1}{\lambda}\}\right)\right), \tag{4.66}$$

the orthogonal projection of $L^2$ onto $\operatorname{Ker}(I-\lambda T)$ from the resolution of the identity for $T$, in the strong operator topology. As a direct consequence of this convergence, it follows using (3.7) and (2.38) that

$$\lim_{n\to\infty}\sup_{s\in\mathbb{R}}\left\|\imath\mu_n\boldsymbol{t}_{|\lambda_n}(s)+\lambda E_T\left(\{\tfrac{1}{\lambda}\}\right)(\boldsymbol{t}(s))\right\| = 0. \tag{4.67}$$

Next, use (3.7), (3.15), (2.5), and (2.9) to write the following chain of relations

$$\begin{aligned}\mu_n\|\widetilde{\boldsymbol{t}}_{m_n|\lambda_n}(s)-\boldsymbol{t}_{|\lambda_n}(s)\| &= \mu_n\|\chi_{m_n}(s)P_{m_n}R^*_{\lambda_n}(\widetilde{T}_{m_n})(\boldsymbol{t}(s))-R^*_{\lambda_n}(T)(\boldsymbol{t}(s))\| \\ &\leqslant \mu_n\left(\widehat{\chi}_{m_n}(s)\|R^*_{\lambda_n}(T)(\boldsymbol{t}(s))\|+\chi_{m_n}(s)\|\left(P_{m_n}R^*_{\lambda_n}(\widetilde{T}_{m_n})-R^*_{\lambda_n}(T)\right)(\boldsymbol{t}(s))\|\right) \\ &= \mu_n\Big(\widehat{\chi}_{m_n}(s)\|R^*_{\lambda_n}(T)(\boldsymbol{t}(s))\|+\|\chi_{m_n}(s)\Big((P_{m_n}-I)R^*_{\lambda_n}(\widetilde{T}_{m_n}) \\ &\quad +(R^*_{\lambda_n}(\widetilde{T}_{m_n})-R^*_{\lambda_n}(T))(I-P_{m_n})+(R^*_{\lambda_n}(\widetilde{T}_{m_n})-R^*_{\lambda_n}(T))P_{m_n}\Big)(\boldsymbol{t}(s))\|\Big) \\ &\leqslant \mu_n\Big(\widehat{\chi}_{m_n}(s)\|R^*_{\lambda_n}(T)(\boldsymbol{t}(s))\|+\chi_{m_n}(s)\|(P_{m_n}-I)R^*_{\lambda_n}(\widetilde{T}_{m_n})(\boldsymbol{t}(s))\| \\ &\quad +\chi_{m_n}(s)\|(R^*_{\lambda_n}(\widetilde{T}_{m_n})-R^*_{\lambda_n}(T))(I-P_{m_n})\,(\boldsymbol{t}(s))\,\| \\ &\quad +|\lambda_n|\,\|R^*_{\lambda_n}(T)(T-\widetilde{T}_{m_n})R^*_{\lambda_n}(\widetilde{T}_{m_n})(\widetilde{\boldsymbol{t}}_{m_n}(s))\|\Big) \\ &\leqslant \mu_n\left(\frac{|\lambda_n|}{\mu_n}\widehat{\chi}_{m_n}(s)\boldsymbol{\tau}(s)+\frac{3\,|\lambda_n|}{\mu_n}\chi_{m_n}(s)\|(P_{m_n}-I)(\boldsymbol{t}(s))\|\right.\end{aligned}$$

$$+\frac{|\lambda_n|^2}{\mu_n}\|(T-\widetilde{T}_{m_n})(\widetilde{\boldsymbol{t}}_{m_n|\lambda_n}(s))\|\Big)$$
$$\leqslant |\lambda_n|\left(\sup_{s\in\widetilde{\mathbb{I}}_{m_n}}\boldsymbol{\tau}(s)+3\sup_{s\in\mathbb{I}_{m_n}}\|(P_{m_n}-I)(\boldsymbol{t}(s))\|\right)+|\lambda_n|^2\varepsilon_n\to 0\quad\text{as } n\to\infty,$$

the convergence to zero in the final line coming from (2.38), (2.22), and (4.65). This convergence together with that in (4.67) imply, via the inequality

$$\|\imath\mu_n\widetilde{\boldsymbol{t}}_{m_n|\lambda_n}(s)+\lambda E_T\left(\left\{\tfrac{1}{\lambda}\right\}\right)(\boldsymbol{t}(s))\|$$
$$\leqslant \mu_n\|\widetilde{\boldsymbol{t}}_{m_n|\lambda_n}(s)-\boldsymbol{t}_{|\lambda_n}(s)\|+\|\imath\mu_n\boldsymbol{t}_{|\lambda_n}(s)+\lambda E_T\left(\left\{\tfrac{1}{\lambda}\right\}\right)(\boldsymbol{t}(s))\|,$$

that

$$\lim_{n\to\infty}\sup_{s\in\mathbb{R}}\|\imath\mu_n\widetilde{\boldsymbol{t}}_{m_n|\lambda_n}(s)+\lambda E_T\left(\left\{\tfrac{1}{\lambda}\right\}\right)(\boldsymbol{t}(s))\|=0.$$

Since (see (4.66) and (4.59)) $\boldsymbol{e}_{|\{\frac{1}{\lambda}\}}(s)=\lambda E_T\left(\left\{\frac{1}{\lambda}\right\}\right)(\boldsymbol{t}(s))$ for all $s\in\mathbb{R}$, the latter shows that the convergence in (4.62) does indeed hold with respect to the $C(\mathbb{R},L^2)$ norm. This convergence, in turn, combined with those displayed in (2.35), (4.61), and (2.33) also implies the desired $C(\mathbb{R}^2,\mathbb{C})$-convergence in (4.63). Indeed, there holds the following:

$$\lim_{n\to\infty}\left\|\imath\mu_n\widetilde{\boldsymbol{T}}_{m_n|\lambda_n}-\boldsymbol{E}_{|\{\frac{1}{\lambda}\}}\right\|_{C(\mathbb{R}^2,\mathbb{C})}$$
$$=\lim_{n\to\infty}\sup_{(s,t)\in\mathbb{R}^2}\Big|\imath\mu_n\lambda_n\langle\widetilde{\boldsymbol{t}}_{m_n}(t),\widetilde{\boldsymbol{t}}_{m_n|\lambda_n}(s)\rangle$$
$$-\lambda^2\left\langle\boldsymbol{t}(t),E_T\left(\left\{\tfrac{1}{\lambda}\right\}\right)(\boldsymbol{t}(s))\right\rangle+\imath\mu_n\widetilde{\boldsymbol{T}}_{m_n}(s,t)\Big|=0,$$

where the first equality uses equation (3.16) and formula (4.60). The theorem is proved. □

For each $n$ in $\mathbb{N}$ and each $\omega$ in $\Omega_0$, let $\boldsymbol{E}_{n|\omega}$ denote the corresponding value of the spectral function for the (also Hermitian) subkernel $\widetilde{\boldsymbol{T}}_n$ of $\boldsymbol{T}$ (see (2.27)), and recall from (4.59) and (4.60) that then, for all $s$, $t\in\mathbb{R}$,

$$\boldsymbol{E}_{n|\omega}(s,t)=\langle\boldsymbol{e}_{n|\omega}(t),\boldsymbol{e}_{n|\omega}(s)\rangle\quad\text{with } \boldsymbol{e}_{n|\omega}(s):=\overline{\boldsymbol{E}_{n|\omega}(s,\cdot)}=v_\omega(\widetilde{T}_n)(\widetilde{\boldsymbol{t}}_n(s)).\tag{4.68}$$

**Theorem 9.** *Suppose that $\omega\in\Omega_0$. Then the sequence $\{\boldsymbol{E}_{n|\omega}\}$ is relatively compact in $C(\mathbb{R}^2,\mathbb{C})$, and hence each sequence of positive integers contains a subsequence $\{n_k\}_{k=1}^{\infty}$ along which $\boldsymbol{E}_{n_k|\omega}$ converges in the norm of $C(\mathbb{R}^2,\mathbb{C})$.*

*Proof.* The proof resembles in outline the proof of Theorem 3. Let $\omega$ belong to $\Omega_0$, and let $M=\sup_{\lambda\in\omega}\frac{1}{|\lambda|}$, so that $\|v_\omega(\widetilde{T}_n)\|\leqslant M$ for all $n$ in $\mathbb{N}$. From Lemma 2 applied with $v_\omega(\widetilde{T}_n)P_n$ in place of $B_n$ and its proof it follows via (4.68) that the estimate

$$\|\boldsymbol{e}_{n|\omega}(x)\|\leqslant M\|\boldsymbol{\tau}\|_{C(\mathbb{R},\mathbb{R})}\tag{4.69}$$

holds for all $x$ in $\mathbb{R}$ and all $n$ in $\mathbb{N}$ (compare (2.41)), and that the subset

$$\mathfrak{e}=\bigcup_{n=1}^{\infty}\boldsymbol{e}_{n|\omega}\subset C(\mathbb{R},L^2)\tag{4.70}$$

is equicontinuous on $\mathbb{R}\cup\{\infty\}$, and, as seen from the proof of (2.42), its equicontinuity is *uniform* on $\mathbb{R}$, in the sense that for a fixed $\varepsilon$, the value of $\delta_x(\varepsilon)$ which is being used in (2.1) can be found to be the same for every $x$ in $\mathbb{R}$. It follows that given any $\varepsilon > 0$, a $\delta > 0$ can be found such that for $x$, $y$, $u$ and $v$ in $\mathbb{R}$,

$$\max\left\{\sup_{n\in\mathbb{N}} \|\boldsymbol{e}_{n|\omega}(x) - \boldsymbol{e}_{n|\omega}(y)\|, \sup_{n\in\mathbb{N}} \|\boldsymbol{e}_{n|\omega}(u) - \boldsymbol{e}_{n|\omega}(v)\|\right\} < \varepsilon$$

whenever $\max\{|x-y|, |u-v|\} < \delta$ or/and whenever both $\min\{|x|, |u|\} > \dfrac{1}{\delta}$ and $\boldsymbol{\tau}(y) = \boldsymbol{\tau}(v) = 0$, whence, for all the $x$, $y$, $u$ and $v$ values so selected, the inequality

$$|\boldsymbol{E}_{n|\omega}(x,u) - \boldsymbol{E}_{n|\omega}(y,v)| \leqslant 2M\|\boldsymbol{\tau}\|_{C(\mathbb{R},\mathbb{R})}\varepsilon \tag{4.71}$$

follows, by means of the more general inequality that makes use of formulae (4.68) and holds for all $x$, $y$, $u$, and $v$ in $\mathbb{R}$ and all $n$ in $\mathbb{N}$:

$$\begin{aligned} &|\boldsymbol{E}_{n|\omega}(x,u) - \boldsymbol{E}_{n|\omega}(y,v)| \\ &\leqslant \|\boldsymbol{e}_{n|\omega}(u)\|\|\boldsymbol{e}_{n|\omega}(x) - \boldsymbol{e}_{n|\omega}(y)\| + \|\boldsymbol{e}_{n|\omega}(y)\|\|\boldsymbol{e}_{n|\omega}(u) - \boldsymbol{e}_{n|\omega}(v)\|. \end{aligned} \tag{4.72}$$

Observe that this inequality, when written for fixed $y$ and $v$ satisfying $\boldsymbol{\tau}(y) = \boldsymbol{\tau}(v) = 0$, also leads (via (4.69)) to the following estimate, valid for all $x$ and $u$ in $\mathbb{R}$ and all $n$ in $\mathbb{N}$:

$$|\boldsymbol{E}_{n|\omega}(x,u)| \leqslant M^2\|\boldsymbol{\tau}\|^2_{C(\mathbb{R},\mathbb{R})}. \tag{4.73}$$

Now, by (4.71), the subset

$$\mathfrak{E} = \bigcup_{n=1}^{\infty} \boldsymbol{E}_{n|\lambda_n} \subset C(\mathbb{R}^2, \mathbb{C}) \tag{4.74}$$

is equicontinuous on $\mathbb{R}^2 \cup \{\infty\}$, i.e., it does have property (A2) of Theorem 1, for $k = 2$ and $B = \mathbb{C}$. By (4.73), the subset $\mathfrak{E}$ in (4.74) does have property (A1) of that theorem, also with $\mathbb{R}^2$ and $\mathbb{C}$ in place of $\mathbb{R}^k$ and $B$, respectively. Linking properties (A1) and (A2) yields, via the same Theorem 1, the relative compactness of $\mathfrak{E}$ in $C(\mathbb{R}^2, \mathbb{C})$. The theorem is proved. □

In the next theorem, we emphasize the case in which the sequences (4.68) converge respectively to desired limits $\boldsymbol{E}_{|\omega}$ and $\boldsymbol{e}_{|\omega}$ without passing to subsequences.

**Theorem 10.** *Let $a$, $b \in \mathbb{R}$ $(a < b)$, and let the open interval $\omega = (a, b)$ be in $\Omega_0$. Then the following limits hold true whenever both $\frac{1}{a} \notin \Lambda_p(T)$ and $\frac{1}{b} \notin \Lambda_p(T)$:*

$$\lim_{n\to\infty} \|\boldsymbol{e}_{n|\omega} - \boldsymbol{e}_{|\omega}\|_{C(\mathbb{R},L^2)} = 0, \tag{4.75}$$

$$\lim_{n\to\infty} \|\boldsymbol{E}_{n|\omega} - \boldsymbol{E}_{|\omega}\|_{C(\mathbb{R}^2,\mathbb{C})} = 0. \tag{4.76}$$

*Hence, by (2.17), $\lim_{n\to\infty} \boldsymbol{E}_{n|\omega}(s,t) = 0$ holds for all $(s,t) \in \mathbb{R}^2$ if and only if $(\frac{1}{b}, \frac{1}{a}) \subset \Pi(T)$ holds.*

*Proof.* For such an interval $\omega = (a, b)$, Theorems VIII.20(b) and VIII.24(b) in [31] ensure that for all $f \in L^2$,

$$v_\omega(\widetilde{T}_n)f = E_{\widetilde{T}_n}(\omega)v^c_\omega(\widetilde{T}_n)f \to E_T(\omega)v^c_\omega(T)f = v_\omega(T)f \quad \text{as } n \to \infty,$$

where $v^c_\omega : \mathbb{R} \to \mathbb{R}$ is a bounded continuous function such that $v^c_\omega(\lambda) = v_\omega(\lambda)$ for every $\lambda \in \omega$. Then, by virtue of (4.59), (4.68), (2.38) and (2.35),

$$\|\boldsymbol{e}_{n|\omega} - \boldsymbol{e}_{|\omega}\|_{C(\mathbb{R},L^2)} \leqslant \sup_{s\in\mathbb{R}} \|(v_\omega(\widetilde{T}_n) - v_\omega(T))(\boldsymbol{t}(s))\| + \|v_\omega(\widetilde{T}_n)\|\|\widetilde{\boldsymbol{t}}_n - \boldsymbol{t}\|_{C(\mathbb{R},L^2)} \to 0$$

as $n \to \infty$, and thus (4.75) holds. Whence, by (4.60) and (4.68),

$$\|\boldsymbol{E}_{n|\omega} - \boldsymbol{E}_{|\omega}\|_{C(\mathbb{R}^2,\mathbb{C})} \leqslant (\|\boldsymbol{e}_{n|\omega}\|_{C(\mathbb{R},L^2)} + \|\boldsymbol{e}_{|\omega}\|_{C(\mathbb{R},L^2)})\|\boldsymbol{e}_{n|\omega} - \boldsymbol{e}_{|\omega}\|_{C(\mathbb{R},L^2)} \to 0$$

as $n \to \infty$, and thus (4.76) also holds. The theorem is proved. □

## References

[1] F. Aràndiga, V. Caselles, "On strongly stable approximations", *Rev. Mat. Univ. Complut. Madrid*, **7**:2 (1994), 207–217.

[2] J. Buescu, "Positive integral operators in unbounded domains", *J. Math. Anal. Appl.*, **296**:1 (2004), 244–255.

[3] T. Carleman, "Zur Theorie der linearen Integralgleichungen", *Math. Z.*, **9** (1921), 196–217.

[4] T. Carleman, *Sur les équations intégrales singulières à noyau réel et symétrique*, A.-B. Lundequistska Bokhandeln, Uppsala, 1923.

[5] F. Chatelin, *Spectral approximation of linear operators*, Academic Press, New York-London, 1983.

[6] C. G. Costley, "On singular normal linear equations", *Can. Math. Bull.*, **13** (1970), 199–203.

[7] N. Dunford, J. T. Schwartz, *Linear operators. Part II: Spectral theory*, Interscience, New York–London, 1963.

[8] N. Dunford, J. T. Schwartz, *Linear operators. Part III: Spectral operators*, Wiley-Interscience, New York, 1971.

[9] S. R. Foguel, "The relations between a spectral operator and its scalar part", *Pacific J. Math.*, **8** (1958), 51–65.

[10] P. Halmos, V. Sunder, *Bounded integral operators on $L^2$ spaces*, Springer, Berlin, 1978.

[11] E. Hille, *Functional analysis and semi-groups*, Amer. Math. Soc. Colloq. Publ. 31, New York, 1948.

[12] K. Jörgens, *Lineare Integraloperatoren*, B. G. Teubner, Stuttgart, 1970 (in German); English transl.: Pitman Advanced Publishing Program, Boston-London-Melbourne, 1982.

[13] L. V. Kantorovich, G. P. Akilov, *Functional analysis*, 3rd ed., Nauka, Moscow, 1984 (in Russian); English transl.: Pergamon Press, Oxford-Elmsford-New York, 1982.

[14] T. Kato, *Perturbation theory for linear operators*, Corr. print. of the 2nd ed., Springer-Verlag, Berlin-Heidelberg-New York, 1980.

[15] V. B. Korotkov, *Integral operators*, Nauka, Novosibirsk, 1983 (in Russian).

[16] V. B. Korotkov, "Some unsolved problems of the theory of integral operators", *Sobolev spaces and related problems of analysis*, Trudy Inst. Mat.. V. 31, Izdat. Ross. Akad. Nauk Sib. Otd. Inst. Mat., Novosibirsk, 1996, 187–196 (in Russian); English transl.: *Siberian Adv. Math.*, **7** (1997), 5–17.

[17] V. B. Korotkov, "On the nonintegrability property of the Fredholm resolvent of some integral operators", *Sibirsk. Mat. Zh.*, **39** (1998), 905–907 (in Russian); English transl.: *Siberian Math. J.*, **39** (1998), 781–783.

[18] V. B. Korotkov, *Introduction to the algebraic theory of integral operators*, Far-Eastern Branch of the Russian Academy of Sciences, Vladivostok, 2000, ISBN: 5-1442-0827-5 (in Russian), 79 pp.

[19] B. Misra, D. Speiser, G. Targonski, "Integral operators in the theory of scattering", *Helv. Phys. Acta*, **36** (1963), 963–980.

[20] S. G. Mikhlin, "On the convergence of Fredholm series", *Doklady AN SSSR*, **XLII**:9 (1944), 374–377 (in Russian).

[21] M. T. Nair, "On strongly stable approximations", *J. Austral. Math. Soc. Ser. A*, **52**:2 (1992), 251–260.

[22] I. M. Novitskii, "Carleman integral equations of the second kind on the half-line", *Mathematical Analysis and Differential Equations*, Interuniv. Collect. Sci. Works, ed. M. M. Lavrent'ev, Novosibirsk State Univ., Novosibirsk, 1987, 187–196 (in Russian).

[23] I. M. Novitskii, "Reduction of linear operators in $L^2$ to integral form with smooth kernels,", *Dokl. Akad. Nauk SSSR*, **318**:5 (1991), 1088–1091 (in Russian); English transl.: *Soviet Math. Dokl.*, **43**:3 (1991), 874–877.

[24] I. M. Novitskii, "Unitary equivalence between linear operators and integral operators with smooth kernels", *Differentsial'nye Uravneniya*, **28**:9 (1992), 1608–1616 (in Russian); English transl.: *Differential Equations*, **28**:9 (1992), 1329–1337.

[25] I. M. Novitskii, "Integral representations of linear operators by smooth Carleman kernels of Mercer type", *Proc. Lond. Math. Soc. (3)*, **68**:1 (1994), 161–177.

[26] I. M. Novitskii, “Fredholm formulae for kernels which are linear with respect to parameter”, *Dal’nevost. Mat. Zh.*, **3**:2 (2002), 173–194 (in Russian).

[27] I. M. Novitskii, “Integral operators with infinitely smooth bi-Carleman kernels of Mercer type”, *Int. Electron. J. Pure Appl. Math.*, **2**:1 (2010), 43–73.

[28] I. M. Novitskii, “Some properties of the resolvent kernels for integral equations with bi-Carleman kernels”, *Dal’nevost. Mat. Zh.*, **16**:2 (2016), 186–208.

[29] A. Pietsch, *Eigenvalues and s-numbers*, Akadem. Verlagsges. Geest & Portig K.-G., Leipzig, 1987.

[30] P. J. Rabier, “Ascoli’s theorem for functions vanishing at infinity and selected applications”, *J. Math. Anal. Appl.*, **290**:1 (2004), 171–189.

[31] M. Reed, B. Simon, *Methods of modern mathematical physics. I. Functional analysis*, rev. ed., Academic Press, San Diego, 1980.

[32] A. F. Ruston, *Fredholm theory in Banach spaces*, Cambridge Tracts in Mathematics, Cambridge Univ. Press, Cambridge e.a., 1986.

[33] F. Smithies, “The Fredholm theory of integral equations”, *Duke Math. J.*, **8** (1941), 107–130.

[34] W. J. Trjitzinsky, “Singular integral equations with complex valued kernels”, *Ann. Mat. Pura Appl.*, **4**:25 (1946), 197–254.

[35] J. W. Williams, “Linear integral equations with singular normal kernels of class I”, *J. Math. Anal. Appl.*, **68**:2 (1979), 567–579.

[36] A. C. Zaanen, “An extension of Mercer’s theorem on continuous kernels of positive type”, *Simon Stevin*, **29** (1952), 113–124.

[37] A. C. Zaanen, *Linear analysis*, North-Holland Publ. Comp., Amsterdam-New York-Oxford, 1964.

*Current address*: Khabarovsk Division, Institute of Applied Mathematics, Far-Eastern Branch of the Russian Academy of Sciences, 54, Dzerzhinskiy Street, Khabarovsk 680 000, RUSSIA

*E-mail address*: `nvigma@mail.ru`

*URL*: `www.iam.khv.ru`